# K-PROCESSES, SCALING LIMIT AND AGING FOR THE TRAP MODEL IN THE COMPLETE GRAPH

By L. R. G. Fontes[1] and P. Mathieu

*Instituto de Mathemática e Estatística and Université de Provence*

We study K-processes, which are Markov processes in a denumerable state space, all of whose elements are stable, with the exception of a single state, starting from which the process enters finite sets of stable states with uniform distribution. We show how these processes arise, in a particular instance, as scaling limits of the trap model in the complete graph, and subsequently derive aging results for those models in this context.

**1. Introduction.** In this paper, we study some properties of a family of Markov processes, which we call K-processes, in particular, their relationship in a special case with the scaling limit of a trap model associated to the Random Energy Model (REM) at low temperature—the trap model in the complete graph, as well as with the aging phenomenon exhibited by that model [7]. These processes are thus prototypes of infinite-volume dynamics for low-temperature (mean-field) spin-glasses.

They have the following remarkable characteristic property. Their state space is countably infinite (we take it to be $\{1, 2, \ldots, \infty\}$), with a single *unstable* state, where by unstable we mean that the process spends zero time at that state at each visit to it; as we will see, that state may be either *instantaneous* or *fictitious* (which are standard terms) in different cases. When in a stable state, the process waits for an exponential time and jumps to the unstable state, starting from which, and here is the striking feature, it enters any finite set of stable states with uniform distribution. In the context of spin-glasses, the stable states represent the low-energy configurations, and the unstable state represents the high-energy configurations. The apparent

Received March 2006; revised July 2007.

[1]Supported in part by Brazilian CNPq Grants 307978/2004-4 and 475833/2003-1, and by FAPESP Grant 04/07276-2.

paradox of the uniformity property is elucidated by a summability condition on the inverse of the jump rates.

It turns out that a class of processes with this uniformity property was introduced by Kolmogorov as an example of a Markov process with an instantaneous state, thus not satisfying his equations in their usual form [21], and it has subsequently been considered by many authors. This class comprises all members of the family we study in the present paper but for an important special case, precisely the one related to the trap model in the complete graph. See Remark 3.2 below for more details.

We have two approaches: an analytical one, based on Dirichlet forms, introduced in Section 2; and one based on an explicit probabilistic construction, in Section 3, at the end of which we argue the equivalence of both points of view. In Section 4, we derive a characterization result for K-processes. Section 5 is devoted to the scaling limit of the trap model in the complete graph, and to deriving an aging result for the associated K-process in this context, which can be seen as an aging result for the trap model in the complete graph itself.

Aging is a dynamical phenomenon observed in disordered systems like spin-glasses at low temperature, signaled by the existence of a limit of a given two-time correlation function of the system started at a high-temperature configuration/state, as both times diverge keeping a fixed ratio between them; the limit should be a nontrivial function of the ratio. This is thus a far-from-equilibrium phenomenon. It has been observed in real spin-glasses and studied extensively in the physics literature. See [8] and references therein.

In [7], a phenomenological model for a Glauber dynamics for the Random Energy Model (REM) is introduced, namely the trap model (in the complete graph), and an aging result for it is established. See more on that model and what is meant by an aging result in Section 5. Roughly speaking, the trap model is a symmetric continuous-time random walk, typically in a regular graph, finite or infinite. The jump rates at the vertices are i.i.d. random variables with a polynomial tail at the origin, whose degree is related to temperature, so that degree less than 1 is equivalent to low temperature. We will assume this regime throughout.

In the mathematics literature, much attention has recently been given to trap models, and many aging results were derived for them. In [4, 5], the trap model in the hypercube is studied, with the rates given by energies of the REM associated to the vertices of the hypercube. The aging result obtained in [5] is for the same correlation function as one considered in [7] with the same limit, thus giving support to the phenomenology underlying the adoption of the trap model in the complete graph by the authors of the latter paper. In [9], an alternative approach to studying the trap model in the hypercube is developed, and the aging result in [7] alluded to above is in this fashion reestablished.

The trap model in $\mathbb{Z}$ was considered in [16] and [1]; the one in $\mathbb{Z}^2$, in [2, 6, 10]; the one in $\mathbb{Z}^d$, $d \geq 3$ in [2, 10]. In [3], a comprehensive approach to obtaining aging results for the trap model on a class of graphs, including $\mathbb{Z}^d$ and tori in two and higher dimensions, the complete graph, the hypercube, is developed.

In most of the above cited work, aging is derived for given correlation functions, without specific regard to the fact that it may arise as a scaling property of the full dynamics. As in [16] and [1], we follow the latter approach for the trap model in the complete graph, and derive its scaling limit (see Theorem 5.2 below); aging results follow (after a further limit is taken, as explained below; see also Theorem 5.11 and Corollary 5.15 below).

It should be noted that, since a time divergence is involved, the scaling limit of the rates (or alternatively the average holding times) should be taken together with the scaling limit of the dynamics, the limiting object acting as a disordered set of parameters for the limiting dynamics. The rescaling is of time only (in such a way that the lowest rates are of order 1), since space is not relevant for the model in the complete graph. The scaling limit results as, roughly speaking, a dynamics in the deepest traps (but the remainder states play a role: they are lumped together in the limit in a single unstable state).

In this model, in the scaling limit, aging is a phenomenon of the dynamics at *vanishing* times: at order 1 or larger times the dynamics is close enough to or in equilibrium, in contrast to the one-dimensional case of [16] and [1], where it could be said that aging occurs for fixed macroscopic times. This should be compared to the aging result in [7] alluded to above, which takes place in a large microscopic time regime (in our case, it occurs at short macroscopic times), and also to the aging result of [3] for the complete graph, taking place at *mesoscopic* time scales; as far as the three regimes can be compared, they coincide, perhaps not surprisingly. See Remarks 5.6 and 5.7.

By taking the scaling limit first, and the aging limit after, we can see aging as a macroscopic phenomenon (taking place in the limiting dynamics). We point out that the latter limit holds for almost every realization of the underlying (macroscopic) disorder: Theorem 5.11 and Corollary 5.15 are almost sure aging results.

The scaling limit for the trap model in the complete graph is not relevant only as a background for aging, even though that is our main motivation for taking it in this paper. It contains also information about other important features of the dynamics at long microscopic times: from aging at short macroscopic times, to approach to equilibrium at large macroscopic times. So it has an interest of its own. Inasmuch as the REM is a prototype for a (mean-field) spin-glass, and the trap model in the complete graph is a prototype for a Glauber dynamics for the REM at low temperature, this scaling limit turns up as a prototype for an infinite-volume dynamics of a (mean-field) spin-glass at low temperature. We expect the same process to arise as an appropriate scaling limit for the trap model in the hypercube (as dimension diverges), and also for the hopping dynamics for the REM, either in the complete graph or the hypercube. It is conceivable that it will also be the scaling limit of the Metropolis dynamics for the REM in the hypercube (see, e.g., [15] for a definition of this dynamics). We also expect variants of the K-process to show up as scaling limits for dynamics of other mean-field models at low temperature, like the GREM, and that they will also exhibit aging.

Our first step in this study is to describe the class of processes that arise as the scaling limit of the trap model in the complete graph. Since they are closely related to the above mentioned class of processes introduced by Kolmogorov through the above mentioned uniformity property, which turns out to characterize the family consisting of both classes (see Section 4 and Theorem 4.1), we chose to start by defining, constructing and studying relevant properties of that larger family, which we refer to as K-processes.

As mentioned above, we do that analytically, through the Dirichlet form associated to the process (in Section 2), and, alternatively, through a probabilistic construction (in Section 3). The former way has the advantage that the K-processes (are reversible and) have quite simple Dirichlet forms, which facilitate the analysis of quantities like the Green function (see Sections 2.2 and 2.3).

The probabilistic construction, besides having its own interest, allows for a direct analysis of the scaling limit for the trap model in the complete graph and the aging issue, without the need of taking transforms (but we do rely on a Tauberian theorem at a specific point of our argument; see proof of Theorem 5.11), and entails the inclusion of more general aging functions in the analysis and results (see Theorem 5.11 and Corollary 5.15), at little extra effort. See Sections 5 and 5.1.

The analytical construction also leads to simple derivations of aging results in a weak sense, after taking Laplace transforms. See the paragraphs starting after Remark 5.20, including Proposition 5.22, before the proof of Theorem 5.11.

In connection with another area of research, as we briefly discuss in Remark 4.2 in Section 4, a K-process can be viewed as a one-point extension of a Markov process beyond its killing time, an object which is of current interest [11, 18].

## 2. Dirichlet forms approach.

2.1. *Construction.* Let $\bar{\mathbb{N}}^*$ be the one-point compactification of $\mathbb{N}^* = \{1, 2, \ldots\}$, with $\infty$ denoting the extra point. In other words, we take $\bar{\mathbb{N}}^*$ with any fixed metric $d$ making it compact. For definiteness, take

$$(2.1) \qquad d(x, y) = |x^{-1} - y^{-1}|, \qquad x, y \in \bar{\mathbb{N}}^*$$

(with $\infty^{-1} = 0$).

Let $\gamma \colon \mathbb{N}^* \to (0, \infty)$ be such that

$$(2.2) \qquad \sum_{x \in \mathbb{N}^*} \gamma(x) < \infty.$$

We extend $\gamma$ to $\bar{\mathbb{N}}^*$ by declaring

$$(2.3) \qquad \gamma(\infty) = 0.$$

Let $\mathcal{C}$ be the space on continuous real-valued functions on $\bar{\mathbb{N}}^*$ and define

$$(2.4) \qquad \mathcal{D} = \left\{ f \colon \bar{\mathbb{N}}^* \to \mathbb{R} \text{ s.t. } \sum_x (f(x) - f(\infty))^2 < \infty \right\}.$$

("$\sum_x$" usually stands for "$\sum_{x \in \mathbb{N}^*}$.") Note that $\mathcal{D}$ is a dense subset of $\mathcal{C}$.

For $f, g \in \mathcal{D}$, consider the bilinear symmetric form

$$(2.5) \qquad \mathcal{E}(f, g) = \sum_x (f(x) - f(\infty))(g(x) - g(\infty)).$$

LEMMA 2.1.   $(\mathcal{E}, \mathcal{D})$ *is a regular Dirichlet form acting on* $L^2(\bar{\mathbb{N}}^*, \gamma)$ *in the sense of* [19].

PROOF.   First note that $\gamma$ has *full support* since we have assumed that $\gamma(x) > 0$ for all $x \in \mathbb{N}^*$.

Clearly $\mathcal{E}$ is bilinear and symmetric. We should check that $\mathcal{D} \subset L^2(\bar{\mathbb{N}}^*, \gamma)$: let $f \in \mathcal{D}$. Without loss of generality, assume that $f(\infty) = 0$. Therefore $\mathcal{E}(f, f) = \sum_x f(x)^2 < \infty$ and $\sum_x f(x)^2 \gamma(x) \le (\sup_x \gamma(x)) \sum_x f(x)^2 < \infty$.

It is easy to check that contractions act on $\mathcal{E}$ so that $\mathcal{E}$ is a Markovian form.

The last point is to prove that $\mathcal{D}$ is *complete* for the norm induced by the bilinear form $\mathcal{E}$: assume that $f_n \in \mathcal{D}$ satisfies $f_n \to 0$ in $L^2(\bar{\mathbb{N}}^*, \gamma)$ and $\mathcal{E}(f_n - f_m, f_n - f_m) \to 0$ as $n$ and $m$ tend to infinity. Then we must have $f_n(x) \to 0$ for any $x \in \mathbb{N}^*$ [because $\gamma(x) > 0$]. Also, for any $\varepsilon > 0$ there exists $n_0$ s.t. for any $n, m \ge n_0$ and any $x \in \mathbb{N}^*$,

$$|f_n(x) - f_n(\infty) - f_m(x) + f_m(\infty)| \le \varepsilon.$$

[This comes from the assumption $\mathcal{E}(f_n - f_m, f_n - f_m) \to 0$.] Letting $m$ go to $\infty$ and then $x$ go to $\infty$, we get that $\limsup_m |f_m(\infty)| \leq \varepsilon$ and therefore $f_m(\infty) \to 0$ as $m$ tends to $\infty$. By Fatou's lemma,

$$\mathcal{E}(f_n, f_n) = \sum_x (f_n(x) - f_n(\infty))^2$$
$$= \sum_x \liminf_m (f_n(x) - f_m(x) - f_n(\infty) + f_m(\infty))^2$$
$$\leq \liminf_m \sum_x (f_n(x) - f_m(x) - f_n(\infty) + f_m(\infty))^2$$
$$= \liminf_m \mathcal{E}(f_n - f_m, f_n - f_m)$$

and therefore $\mathcal{E}(f_n, f_n) \to 0$.  $\square$

REMARK 2.2. We first recall from Chapter 1 of [19] that, to any Dirichlet form, can be associated a Markovian semigroup. Thus there exists a strongly continuous semi-group of symmetric, Markovian contractions of $L^2(\bar{\mathbb{N}}^*, \gamma)$, say $(P_t, t \geq 0)$, whose Dirichlet form is $(\mathcal{E}, \mathcal{D})$ in the sense that, for any $f, g \in \mathcal{D}$,

$$\frac{1}{t} \sum_x (f(x) - P_t f(x)) g(x) \gamma(x) \to \mathcal{E}(f, g).$$

Since $(\mathcal{E}, \mathcal{D})$ is also regular, we may apply Theorem 7.2.1 of [19] to conclude that there also exists a symmetric Markov process, in fact a Hunt process, whose Dirichlet form is $(\mathcal{E}, \mathcal{D})$ on $L^2(\bar{\mathbb{N}}^*, \gamma)$: there exists a Markovian family of probability measures on the space of càdlàg trajectories in $\bar{\mathbb{N}}^*$, say $(\mathbb{P}_x, x \in \bar{\mathbb{N}}^*)$, with the property that

$$(2.6) \qquad\qquad \mathbb{E}_x(f(X_t)) = P_t f(x),$$

for $t \geq 0$, any $x \in \bar{\mathbb{N}}^*$ and any $f \in L^2(\bar{\mathbb{N}}^*, \gamma)$. (In the above expression, $X_t$ is the coordinate map and $\mathbb{E}_x$ is the expectation with respect to $\mathbb{P}_x$.) We refer to Appendix A.2 of [19] for an introduction to Hunt processes. The Markov process $(\mathbb{P}_x, x \in \bar{\mathbb{N}}^*)$ is uniquely determined by the Dirichlet form $(\mathcal{E}, \mathcal{D})$ in the sense of Theorem 4.2.7 of [19].

2.2. *Computation of hitting times and capacities.* Given the explicit enough form of $\mathcal{E}$ it is easy to compute the law of some hitting times and entrance laws.

Let $A$ be a subset of $\bar{\mathbb{N}}^*$ and $\tau^A = \inf\{t\,; X(t) \in A\}$. Let $\mathcal{F}_A$ be the set of functions in $\mathcal{D}$ that vanish on $A$. Let $f : \bar{\mathbb{N}}^* \to \mathbb{R}$ be bounded measurable and choose $\lambda > 0$. Then the function $y \to \mathbb{E}_y(f(X_{\tau^A})e^{-\lambda \tau^A})$ is the orthogonal projection of the function $f$ on the orthogonal complement of $\mathcal{F}_A$ when $\mathcal{D}$ is equipped with the Hilbert norm $\mathcal{E}(u, u) + \lambda \gamma(u^2)$; see Theorem 4.3.1 of [19]. We will repeatedly use this fact to make explicit computations in the next lemmas

LEMMA 2.3.  *Let $\tau_x = \inf\{t \,;\, X(t) \neq x\}$. Then*

$$(2.7) \qquad \mathbb{E}_x(e^{-\lambda \tau_x}) = \frac{1}{1 + \lambda \gamma(x)}.$$

PROOF.  The function $y \to \mathbb{E}_y(e^{-\lambda \tau_x})$ is the minimizer of the expression $\mathcal{E}(u, u) + \lambda \gamma(u^2)$ among functions $u$ satisfying $u(y) = 1$ for $y \neq x$. But, for such a function $u$, we have $\mathcal{E}(u, u) + \lambda \gamma(u^2) = (u(x) - 1)^2 + \lambda \gamma(x) u(x)^2 + \lambda(1 - \gamma(x))$ that is minimal for $u(x) = \frac{1}{1 + \lambda \gamma(x)}$.  □

LEMMA 2.4.  *Let $\sigma_\infty = \inf\{t \,;\, X(t) = \infty\}$. Then*

$$(2.8) \qquad \mathbb{E}_x(e^{-\lambda \sigma_\infty}) = \frac{1}{1 + \lambda \gamma(x)}.$$

PROOF.  We now have to minimize $\mathcal{E}(u, u) + \lambda \gamma(u^2)$ among functions $u$ satisfying $u(\infty) = 1$. But for such a function $u$, we have $\mathcal{E}(u, u) + \lambda \gamma(u^2) = \sum_x (u(x) - 1)^2 + \lambda \gamma(x) u(x)^2$ that is minimal for $u(x) = \frac{1}{1 + \lambda \gamma(x)}$.  □

REMARK 2.5.  In particular note that $\sigma_\infty < \infty$ $\mathbb{P}_x$ a.s. Hence $\mathbb{P}_\infty$ is well defined. Since $\mathbb{E}_x(e^{-\lambda \sigma_\infty}) = \mathbb{E}_x(e^{-\lambda \tau_x})$ and $\tau_x \leq \sigma_\infty$, we must have $\tau_x = \sigma_\infty$ $\mathbb{P}_x$ a.s. In particular $X(\tau_x) = \infty$ $\mathbb{P}_x$ a.s.

LEMMA 2.6.  *Let $A$ be a finite subset of $\mathbb{N}^*$ of size $n$, and $\tau^A = \inf\{t \,;\, X(t) \in A\}$. Then, for any function $f : A \to \mathbb{R}$, any $\lambda > 0$ and any $y \notin A$, we have*

$$(2.9) \begin{aligned} &\mathbb{E}_y(f(X_{\tau^A}) e^{-\lambda \tau^A}) \\ &\quad = \left(n + (1 + \lambda \gamma(y)) \sum_{x \notin A} \lambda \gamma(x)/(1 + \lambda \gamma(x))\right)^{-1} \sum_{x \in A} f(x). \end{aligned}$$

*In particular, for $\lambda = 0$, we find that the law of $X(\tau^A)$ is uniform over $A$.*

PROOF.  We have to minimize $\mathcal{E}(u, u) + \lambda \gamma(u^2)$ among functions $u$ satisfying $u(x) = f(x)$ for $x \in A$. For such a function $\mathcal{E}(u, u) + \lambda \gamma(u^2) = \sum_{x \in A} (f(x) - u(\infty))^2 + \sum_{x \notin A} (u(x) - u(\infty))^2 + \lambda \gamma(x) u(x)^2 + \lambda \sum_{x \in A} \gamma(x) f(x)^2$. The solution has the form $u(y) = \frac{u(\infty)}{1 + \lambda \gamma(y)}$ for $y \notin A$ and we find $u(\infty)$ by minimizing $\sum_{x \in A} (f(x) - u(\infty))^2 + u(\infty)^2 \sum_{x \notin A} \frac{\lambda \gamma(x)}{1 + \lambda \gamma(x)}$.  □

After a similar computation, we get the following.

LEMMA 2.7. *Let $A$ be as in the previous lemma. Then*

$$(2.10) \quad \mathbb{E}_\infty(f(X_{\tau^A})e^{-\lambda\tau^A}) = \frac{1}{n + \sum_{x\notin A}\lambda\gamma(x)/(1+\lambda\gamma(x))}\sum_{x\in A}f(x).$$

It is also possible to compute the Green kernel

$$(2.11) \quad g_\lambda(x) = \lambda\int_0^\infty e^{-\lambda s}\mathbb{P}_\infty(X(s)=x)\,ds.$$

The Markov property gives

$$(2.12) \quad g_\lambda(x) = \mathbb{E}_\infty(e^{-\lambda\tau^{\{x\}}})(1-\mathbb{E}_x(e^{-\lambda\tau_x})) + \mathbb{E}_\infty(e^{-\lambda\tau^{\{x\}}})\mathbb{E}_x(e^{-\lambda\tau_x})g_\lambda(x).$$

[Remember that $X(\tau_x)=\infty$ a.s.] Using Lemma 2.6, we get that

$$(2.13) \quad g_\lambda(x) = \frac{\lambda\gamma(x)/(1+\lambda\gamma(x))}{\sum_y\lambda\gamma(y)/(1+\lambda\gamma(y))}.$$

We also have the following more general formula. Let $g_\lambda(x,y) = \lambda\int_0^\infty e^{-\lambda s}\times\mathbb{P}_y(X(s)=x)\,ds$. Then

$$(2.14) \quad g_\lambda(x,y) = \frac{1}{1+\lambda\gamma(y)}g_\lambda(x).$$

The last formula describes some correlation function whose definition is motivated by so-called aging.

LEMMA 2.8. *Let*

$$(2.15) \quad c_\lambda(\mu) = \int_0^\infty \lambda e^{-\lambda s}\,ds\int_0^\infty \mu e^{-\mu t}\,dt\,\mathbb{P}_\infty(X(u)=X(s)\,\forall u\in[s,s+t]).$$

*Then*

$$(2.16) \quad c_\lambda(\mu) = \frac{\sum_x(\lambda\gamma(x)/(1+\lambda\gamma(x)))(\mu\gamma(x)/(1+\mu\gamma(x)))}{\sum_x\lambda\gamma(x)/(1+\lambda\gamma(x))}.$$

PROOF. As for the Green function, we use the Markov property to write that

$$\begin{aligned}
c_\lambda(\mu) &= \sum_x\int_0^\infty \lambda e^{-\lambda s}\,ds\int_0^\infty \mu e^{-\mu t}\,dt\,\mathbb{P}_\infty(X(s)=x;X(u)=x\,\forall u\in[s,s+t])\\
&= \sum_x\int_0^\infty \lambda e^{-\lambda s}\,ds\int_0^\infty \mu e^{-\mu t}\,dt\,\mathbb{P}_\infty(X(s)=x)\mathbb{P}_x(\tau_x>t)\\
&= \sum_x g_\lambda(x)(1-\mathbb{E}_x(e^{-\mu\tau_x})). \qquad\qquad\qquad\qquad\qquad\square
\end{aligned}$$

2.3. *Some extension.* Let $c > 0$ and define the new measure $\gamma^c = \gamma + c\delta_\infty$. The bilinear form $(\mathcal{E}, \mathcal{D})$ turns out to be also a Dirichlet form when acting on $L^2(\bar{\mathbb{N}}^*, \gamma^c)$. The corresponding Markov process can be described as follows: let $L(t)$ be the local time of $X$ at $\infty$. [$L(t)$ is the unique additive functional whose Revuz measure is $\delta_\infty$.] Define

$$(2.17) \qquad A^c(t) = t + cL(t) \quad \text{and} \quad X^c(t) = X(A^{-1}(t)).$$

Then, under $\mathbb{P}_x$, $X^c$ is a Markov process and its Dirichlet form is $(\mathcal{E}, \mathcal{D})$ acting on $L^2(\bar{\mathbb{N}}^*, \gamma^c)$. Call $\mathbb{P}_x^c$ its law when starting from $x$.

One can then reproduce the same computation as before. In particular we get the expression of the Green function:

$$(2.18) \qquad g_\lambda^c(x) = \frac{\lambda\gamma(x)/(1 + \lambda\gamma(x))}{c\lambda + \sum_y \lambda\gamma(y)/(1 + \lambda\gamma(y))}, \qquad x \in \mathbb{N}^*$$

and, since $g_\lambda^c(\infty) = 1 - \sum_x g_\lambda^c(x)$,

$$(2.19) \qquad g_\lambda^c(\infty) = \frac{c\lambda}{c\lambda + \sum_x \lambda\gamma(x)/(1 + \lambda\gamma(x))}.$$

We finally have that

$$(2.20) \qquad g_\lambda^c(x, y) = \frac{1}{1 + \lambda\gamma(y)} g_\lambda^c(x), \qquad x, y \in \bar{\mathbb{N}}^*.$$

REMARK 2.9. It follows from the character of the time change (2.17) that for all $c \geq 0$, at the entrance time of finite subsets $A$ of $\mathbb{N}^*$ by $X^c$, starting outside $A$, its distribution is uniform in $A$. See [19], Section 6.2, in particular Theorem 6.2.1.

**3. Probabilistic point of view.** In this section we make an explicit construction for the processes introduced in the previous section, and study some of its properties which are relevant for what follows.

Let $\mathcal{N} = \{(N_t^{(x)})_{t \geq 0}, x \in \mathbb{N}^*\}$ be i.i.d. Poisson processes of rate 1, with $\sigma_j^{(x)}$ the $j$th event time of $N^{(x)}$, and $\mathcal{T} = \{T_0; T_i^{(x)} i \geq 1, x \in \mathbb{N}^*\}$ be i.i.d. exponential random variables of rate 1. $\mathcal{N}$ and $\mathcal{T}$ are assumed independent.

For $c \geq 0$ and $y \in \bar{\mathbb{N}}^*$, let

$$(3.1) \qquad \Gamma(t) = \Gamma^{c,y}(t) = \gamma(y) T_0 + \sum_{x=1}^\infty \gamma(x) \sum_{i=1}^{N_t^{(x)}} T_i^{(x)} + ct,$$

where, by convention, $\sum_{i=1}^0 T_i^{(x)} = 0$ for every $x$.

Let $c \geq 0$ be fixed. We define the process $\tilde{X}^{c,y}$ on $\bar{\mathbb{N}}^*$ starting at $y \in \bar{\mathbb{N}}^*$ as follows. For $t \geq 0$

$$(3.2) \qquad \tilde{X}^c(t) = \tilde{X}^{c,y}(t) = \begin{cases} y, & \text{if } 0 \leq t < \gamma(y)T_0, \\ x, & \text{if } \Gamma(\sigma_j^{(x)}-) \leq t < \Gamma(\sigma_j^{(x)}) \\ & \text{for some } 1 \leq j < \infty, \\ \infty, & \text{otherwise.} \end{cases}$$

DEFINITION 3.1. We call $\tilde{X}^{c,y}$ the K-*process* with parameters $\gamma$ and $c$. We will also call it sometimes the K($\gamma$, $c$)-process for shortness.

REMARK 3.2. The case $c = 1$ was introduced by Kolmogorov [21] as an example of a Markov process in a countable state space with an instantaneous state. It is known in this context as the *first example of Kolmogorov* or K1 (Kolmogorov also introduced a second such example, known as K2, which is *not* a K-process by our definition for any $c \geq 0$ and $\gamma$). The case $c = 1$ was then studied in [20] and [12] (Example 3 in Part II, Chapter 20 of the latter reference), where an equivalent construction to the above one is given, and elsewhere (e.g., [17]). The general case of $c > 0$ is not really different from the one introduced by Kolmogorov; one can go from one case to the other by a uniform deterministic time rescaling. The $c = 0$ case is already considerably different. For one thing, it is not strongly continuous (a Markov process $Y$ in $\bar{\mathbb{N}}^*$ is said to be strongly continuous if $\lim_{t\to 0} \mathbb{P}_x(Y_t = y) = \delta_{xy}$, the Kronecker's delta, for all $x, y \in \bar{\mathbb{N}}^*$; see [13], Chapter 2. As result of Lemma 3.15 below, this property is seen to fail for the K-process with $c = 0$ for $x = y = \infty$), which the K1 process is; following Lévy's classification [22], the K-process is of the fourth kind for $c = 0$, and of the fifth kind for $c > 0$. (A process in a countable state space with one unstable state is of the fourth kind if the set of times when the process visits that state is an uncountable set of null Lebesgue measure, and of the fifth kind if that set is a Cantor set of positive Lebesgue measure; see [22], Chapter II. According to this classification, $\infty$ is termed a *fictitious* state when $c = 0$.) Even though the $c = 0$ case is a natural extension of the $c > 0$ one, we did not find any explicit mention of it in the literature. (In [22], though, it is argued in general terms that by looking at a fifth-kind process outside the instantaneous state, one gets a fourth-kind process.) Nevertheless we will show that precisely the $c = 0$ case arises as the scaling limit of a (mean-field) disordered spin dynamics (the trap model in the complete graph) at low temperatures. Its irregular behavior near $\infty$, associated in particular to its lack of strong continuity, is behind the *aging phenomenon* exhibited by such dynamics at such temperatures [7] (see Remark 5.14 below).

REMARK 3.3.   It is clear that $\tilde{X}^{c,y}(0) = y$ almost surely for all $y \in \mathbb{N}^*$. That this also holds for $y = \infty$ follows readily from (3.2).

REMARK 3.4.   Note on the one hand that $T_0$, $\Gamma(\sigma_j^{(x)}-)$, $\Gamma(\sigma_j^{(x)})$ are continuous random variables for every $x \in \bar{\mathbb{N}}^*$ and $j \geq 1$, and on the other hand that $\tilde{X}^{c,y}$ is almost surely continuous off $\{\gamma(y)\,T_0; \Gamma(\sigma_j^{(x)}-), \Gamma(\sigma_j^{(x)}), x \in \bar{\mathbb{N}}^*, j \geq 1\}$. These readily imply that every $s \geq 0$ is almost surely a continuity point of $\tilde{X}^{c,y}$.

REMARK 3.5.   It readily follows from (3.2) that

$$(3.3) \qquad \tilde{X}^{c,y}(t) = \tilde{X}^{c,\infty}(t - \gamma_y\,T_0) \qquad \text{for } t \geq \gamma_y\,T_0.$$

PROPOSITION 3.6.   $\tilde{X}^c$ is càdlàg and Markovian.

REMARK 3.7.   A treatment of the case $c = 1$ can be found in [20] and [12]. Even though both have a construction equivalent to ours, complete proofs of some key properties of the constructed process, like the Markov one, are not presented. For this reason, and in order to include the $c = 0$ case as well, we present below a proof of Proposition 3.6.

The proof is based on strongly approximating $\tilde{X}^c$ in Skorohod space by Markov processes that we now define. For $n \geq 1$ and $y \in \{1, \ldots, n, \infty\}$, let

$$(3.4) \qquad \Gamma_n(t) = \Gamma_n^{c,y}(t) = \gamma(y)\,T_0 + \sum_{x=1}^{n} \gamma(x) \sum_{i=1}^{N_t^{(x)}} T_i^{(x)} + ct$$

and

$$(3.5) \qquad \tilde{X}_n^{c,y}(t) = \begin{cases} y, & \text{if } 0 \leq t < \gamma(y)\,T_0, \\ x, & \text{if } \Gamma_n(\sigma_j^{(x)}-) \leq t < \Gamma_n(\sigma_j^{(x)}) \\ & \text{for some } 1 \leq x \leq n,\ j \geq 1, \\ \infty, & \text{otherwise.} \end{cases}$$

REMARK 3.8.   We note that $\tilde{X}_n^{0,y}$ never visits $\infty$, even when $y = \infty$. See next remark.

REMARK 3.9.   The order in which the sites of $\{1, \ldots, n\}$ are visited by $\tilde{X}_n^{c,y}$ (in case $y$ is finite, after leaving the initial state) is given by the respective (chronological) order of $\{\sigma_j^{(x)}; 1 \leq x \leq n, j \geq 1\}$. Let us denote the latter set by $\mathcal{S}^n = \{S_1^n, S_2^n, \ldots\}$, with $S_1^n < S_2^n < \cdots$. Then $\mathcal{S}^n$ is a Poisson point process of rate $n$, each point of which is labeled according to a different element of an i.i.d. family of uniform in $\{1, \ldots, n\}$ random variables. This implies that the jump probabilities of $\tilde{X}_n^{c,y}$ from any site in case $c = 0$, and from $\infty$ in case $c > 0$, are uniform in $\{1, \ldots, n\}$, and also implies that $\tilde{X}_n^{0,\infty}(0)$ is uniformly distributed in $\{1, \ldots, n\}$ (since it is the label of $S_1^n$; see previous remark).

In case $c > 0$, $c(S_i^n - S_{i-1}^n)$, $i \geq 1$, where $S_0^n \equiv 0$, represent the successive holding times at $\infty$. It is clear then that these times form an i.i.d. sequence of exponential random variables of mean $c/n$.

We have the following two results.

LEMMA 3.10.   $\tilde{X}_n^c$ is càdlàg and Markovian for every $n \geq 1$ and $y \in \{1, \ldots, n, \infty\}$.

LEMMA 3.11.   $\tilde{X}_n^c \to \tilde{X}^c$ as $n \to \infty$ almost surely in the Skorohod norm for every $y \in \bar{\mathbb{N}}^*$.

PROOF OF THE FIRST ASSERTION OF PROPOSITION 3.6.   The first assertions of Lemmas 3.10 and 3.11 readily establish the first assertion of Proposition 3.6 (see [13]).   □

PROOF OF LEMMA 3.10.   Let the starting point $y$ be fixed.

For $c = 0$, $\tilde{X}_n^c$ is the following Markov process on $\{1, \ldots, n\}$. $\tilde{X}_n^{c,y}$ starts at $y$ if $y \in \{1, \ldots, n\}$; $\tilde{X}_n^{c,\infty}$ has uniform initial distribution. When at $x \in \{1, \ldots, n\}$, it waits an exponential time of mean $\gamma(x)$ and then jumps uniformly at random to a site in $\{1, \ldots, n\}$ (which could be $x$ again). See Remarks 3.8 and 3.9 above.

For $c > 0$, $\tilde{X}_n^c$ is the following Markov process on $\{1, \ldots, n, \infty\}$. $\tilde{X}_n^{c,y}$ starts at $y$. When at $x \in \{1, \ldots, n\}$, it waits an exponential time of mean $\gamma(x)$ and then jumps deterministically to $\infty$. When at $\infty$, it waits an exponential time of mean $c/n$, and then jumps uniformly at random to a site in $\{1, \ldots, n\}$. See Remark 3.9 above.   □

PROOF OF LEMMA 3.11.   Let $y$ be fixed, and suppose $n \geq y$ if $y \in \mathbb{N}^*$. We show the almost sure validity of (c) of Proposition 5.3 in Chapter 3 of [13] (page 119).

For $m \in \mathbb{N}^*$, let $\delta_m = \mathrm{diam}\{x \in \bar{\mathbb{N}}^* : x > m\} = (m+1)^{-1}$ and $\{S_1^m < S_2^m < \cdots\} = \{\sigma_j^{(x)}, j \geq 1, 1 \leq x \leq m\}$, with the latter being well defined almost surely. Fix $T > 0$ and let

$$(3.6) \qquad L_n^m = \min\{i \geq 1 : \Gamma_n(S_i^m) \geq T\},$$

which is almost surely finite, and make $S_0^m \equiv 0$. Notice that $L_n^m = L_n^m(y)$ is nonincreasing in $\gamma(y)$ where $y$ is the starting point, and thus

$$(3.7) \qquad \max_{y \in \mathbb{N}^*} L_n^m(y) = L_n^m(\infty).$$

We can now almost surely find $n_m$ so large that $\min_{0 \le i \le L_n^m - 1} [\Gamma_n(S_{i+1}^m -) - \Gamma_n(S_i^m)] > 0$ for $n \ge n_m$. (We can take $n_m \equiv 1$ when $c > 0$.) For these $n$ then define $\lambda_n^m \colon [0, \Gamma_n(S_{L_n^m}^m)] \to \mathbb{R}^+$ inductively as follows:

$$(3.8) \qquad \lambda_n^m(t) = t \qquad \text{if } 0 \le t < \gamma(y) T_0,$$

and for $0 \le i \le L_n^m - 1$ and $\Gamma_n(S_i^m) \le t \le \Gamma_n(S_{i+1}^m)$, let

$$(3.9) \qquad \lambda_n^m(t) = \begin{cases} \Gamma(S_i^m) + \dfrac{\Gamma(S_{i+1}^m -) - \Gamma(S_i^m)}{\Gamma_n(S_{i+1}^m -) - \Gamma_n(S_i^m)} [t - \Gamma_n(S_i^m)], \\ \qquad\qquad\qquad\qquad\qquad \text{if } \Gamma_n(S_i^m) \le t \le \Gamma_n(S_{i+1}^m -), \\ \Gamma(S_{i+1}^m -) - \Gamma_n(S_{i+1}^m -) + t, \quad \text{if } \Gamma_n(S_{i+1}^m -) \le t \le \Gamma_n(S_{i+1}^m). \end{cases}$$

It has the following properties. For all $T > 0$, $m \in \mathbb{N}^*$ and $n \ge m \vee n_m \vee y$

$$(3.10) \qquad \lambda_n^m(t) \ge t, \qquad 0 \le t \le T,$$

$$(3.11) \quad \sup_{0 \le t \le T} |\lambda_n^m(t) - t| \le \max\{\Gamma(S_{i+1}^m -) - \Gamma_n(S_{i+1}^m -); 0 \le i \le L_n^m(\infty) - 1\}$$

[where we have made use of (3.7)], and

$(3.12)$ the right-hand side of (3.11) vanishes almost surely as $n \to \infty$.

Furthermore,

$$(3.13) \qquad \sup_{0 \le t \le T} \text{dist}(\tilde{X}^c(\lambda_n^m(t)), \tilde{X}_n^c(t)) \le \delta_m,$$

since for $t \in [0, T]$, $\tilde{X}^c(\lambda_n^m(t))$ and $\tilde{X}_n^c(t)$ coincide when either one is in $\{1, \ldots, m\}$.

From (3.11), for every $m \in \mathbb{N}^*$ there almost surely exists $n_m' \ge n_m$ such that for $n \ge n_m'$

$$(3.14) \qquad \sup_{0 \le t \le T} |\lambda_n^m(t) - t| \le \delta_m$$

and (3.13) hold. We may assume $(n_m')$ is strictly increasing.

For $n \ge n_1'$, let $m_n = i$ when $n_i' \le n < n_{i+1}'$ and $\tilde{\lambda}_n = \lambda_n^{m_n}$. We then have for $T > 0$

$$(3.15) \qquad \sup_{0 \le t \le T} \text{dist}(\tilde{X}^c(\tilde{\lambda}_n(t)), \tilde{X}_n^c(t)) \to 0,$$

$$(3.16) \qquad \sup_{0 \le t \le T} |\tilde{\lambda}_n(t) - t| \to 0,$$

almost surely as $n \to \infty$, and the above mentioned condition (c) is verified. $\square$

REMARK 3.12. Since the right-hand side of (3.11) is independent of $y \in \bar{\mathbb{N}}^*$, the convergence in (3.16) is actually uniform in $y \in \bar{\mathbb{N}}^*$.

We will also need the following lemma to prove the second assertion of Proposition 3.6.

LEMMA 3.13. *For every $t \geq 0$, $\tilde{X}^{c,y}(t) \to \tilde{X}^{c,\infty}(t)$ as $y \to \infty$ almost surely.*

PROOF. The case $t = 0$ is clear. For $t > 0$, since we are taking $y \to \infty$, we may assume that $\gamma(y) T_0 \leq t$ and then from (3.3) we have that

$$(3.17) \qquad |\tilde{X}^{c,y}(t) - \tilde{X}^{c,\infty}(t)| = |\tilde{X}^{c,\infty}(t - \gamma(y) T_0) - \tilde{X}^{c,\infty}(t)|$$

and the result follows from Remark 3.4.  □

PROOF OF THE SECOND ASSERTION OF PROPOSITION 3.6. We will show that the Markov property of $\tilde{X}_n^c$ survives in the limit. Before going into the argument, it is perhaps worth mentioning that not all limits of Markov processes are Markovian. In the present case, a few regularity properties enjoyed by the processes involved are behind the fact (or, to be more precise, behind our argument below), namely the Feller property of $\tilde{X}^c$ (see Remark 3.14 below), the almost sure continuity of $\tilde{X}^c(s)$ at all deterministic $s$ (see Remark 3.4 above) and the uniformity in $y \in \bar{\mathbb{N}}^*$ of the convergence of $\tilde{\lambda}_n(t)$ to $t$ (see Remark 3.12 above). The argument below is perhaps a bit technical; all the above mentioned properties enter it in due course.

Now for the argument. Lemma 3.10 implies that for arbitrary $m \geq 1$, $0 \leq t_1 < \cdots < t_{m+1}$ and bounded continuous functions $f_1, \ldots, f_{m+1}$, we have that

$$(3.18) \qquad \begin{aligned} &\mathbb{E}[f_1(\tilde{X}_n^c(t_1)) \cdots f_m(\tilde{X}_n^c(t_m)) f_{m+1}(\tilde{X}_n^c(t_{m+1}))] \\ &\qquad = \mathbb{E}[f_1(\tilde{X}_n^c(t_1)) \cdots f_m(\tilde{X}_n^c(t_m)) \Psi_{t_{m+1}-t_m}^n f_{m+1}(\tilde{X}_n^c(t_m))], \end{aligned}$$

where $\Psi^n$ is the semigroup of $\tilde{X}_n^c$, that is, for $t \geq 0$, a continuous function $f$ and $y \in \bar{\mathbb{N}}^*$

$$(3.19) \qquad\qquad \Psi_t^n f(y) = \mathbb{E}[f(\tilde{X}_n^{c,y}(t))].$$

By Lemma 3.11, the left-hand side of (3.18) converges to

$$(3.20) \qquad \mathbb{E}[f_1(\tilde{X}^c(t_1)) \cdots f_m(\tilde{X}^c(t_m)) f_{m+1}(\tilde{X}^c(t_{m+1}))]$$

as $n \to \infty$. Let us estimate the right-hand side of (3.18) by

$$(3.21) \quad \mathbb{E}[f_1(\tilde{X}_n^c(t_1)) \cdots f_m(\tilde{X}_n^c(t_m)) \Psi_{t_{m+1}-t_m} f_{m+1}(\tilde{X}_n^c(t_m))] + \epsilon_n,$$

where for $t \geq 0$, a continuous function $f$ and $y \in \bar{\mathbb{N}}^*$

$$(3.22) \qquad \Psi_t f(y) = \mathbb{E}[f(\tilde{X}^{c,y}(t))]$$

and

$$(3.23) \qquad |\epsilon_n| \leq \text{const} \sup_y |\Psi^n_{t_{m+1}-t_m} f_{m+1}(y) - \Psi_{t_{m+1}-t_m} f_{m+1}(y)|.$$

From Lemma 3.13, we have that $\Psi_{t_{m+1}-t_m} f_{m+1}(\cdot)$ is continuous, and now Lemma 3.11 implies that the left term of (3.21) converges to

$$(3.24) \qquad \mathbb{E}[f_1(\tilde{X}^c(t_1)) \cdots f_m(\tilde{X}^c(t_m)) \Psi_{t_{m+1}-t_m} f_{m+1}(\tilde{X}^c(t_m))]$$

as $n \to \infty$.

Let us now examine the right-hand side of (3.23). We first relabel $f_{m+1} = g$, $t_{m+1} - t_m = s$ and $\gamma(y) = \gamma_y$. We have that

$$(3.25) \qquad \begin{aligned} &\Psi^n_{t_{m+1}-t_m} f_{m+1}(y) - \Psi_{t_{m+1}-t_m} f_{m+1}(y) \\ &= \mathbb{E}[g(\tilde{X}^{c,y}_n(s)) - g(\tilde{X}^{c,y}(s))] \\ &= \mathbb{E}[g(\tilde{X}^{c,y}_n(s)) - g(\tilde{X}^{c,y}(\tilde{\lambda}_n(s)))] \\ &\quad + \mathbb{E}[g(\tilde{X}^{c,y}(\tilde{\lambda}_n(s))) - g(\tilde{X}^{c,y}(s))], \end{aligned}$$

with $\tilde{\lambda}_n$ as defined in the paragraph of (3.15), with $T > s$. From (3.13), it follows that the sup in $y$ of the absolute value of the first expected value in the right-hand side of (3.25) vanishes as $n \to \infty$ (since $g$ is continuous, and thus uniformly continuous since $\bar{\mathbb{N}}^*$ is compact).

Remark 3.4 now implies that there exists a sequence $k_n$ going to infinity as $n \to \infty$ such that as $n \to \infty$

$$(3.26) \qquad \max_{1 \leq y \leq k_n} |\mathbb{E}[g(\tilde{X}^{c,y}(\tilde{\lambda}_n(s))) - g(\tilde{X}^{c,y}(s))]| \to 0.$$

We now note that, from (3.10), (3.12), $\tilde{\lambda}_n(s) \geq s$, and that $\tilde{\lambda}_n(s) \to s$ as $n \to \infty$ uniformly in $y$ almost surely (see Remark 3.12 above). From this and (3.3) we then have

$$(3.27) \qquad \begin{aligned} &\sup_{y > k_n} |\mathbb{E}[g(\tilde{X}^{c,y}(\tilde{\lambda}_n(s))) - g(\tilde{X}^{c,y}(s))]| \\ &\leq \sup_{y > k_n} \mathbb{E}[|g(\tilde{X}^{c,\infty}(\tilde{\lambda}_n(s) - \gamma_y T_0)) \\ &\qquad\qquad - g(\tilde{X}^{c,\infty}(s - \gamma_y T_0))|; \gamma_y T_0 < \delta_n] \\ &\quad + 2\|g\| \sup_{y > k_n} \mathbb{P}(\gamma_y T_0 \geq \delta_n), \end{aligned}$$

where $\delta_n$ is chosen so that $\delta_n \to 0$ and $\sup_{y > k_n} \mathbb{P}(\gamma_y T_0 \geq \delta_n) = e^{-\delta_n / (\sup_{y > k_n} \gamma_y)} \to 0$ as $n \to \infty$. Then we have that the latter summand in the right-hand side of (3.27) vanishes as $n \to \infty$. The former one can be bounded above by

$$\mathbb{E}\left[\sup_{y > k_n} \sup_{0 \leq t \leq \delta_n} |g(\tilde{X}^{c,\infty}(\tilde{\lambda}_n(s) - t)) - g(\tilde{X}^{c,\infty}(s - t))|\right]$$

and as $n \to \infty$ that vanishes as well by Remark 3.4 above [(any) $s$ is almost surely a continuity point of $\tilde{X}^{c,\infty}(\cdot)$].

We have thus concluded that $|\epsilon_n| \to 0$ as $n \to \infty$, and then from (3.18)–(3.24) we have that

$$(3.28) \quad \begin{aligned} &\mathbb{E}[f_1(\tilde{X}^c(t_1)) \cdots f_m(\tilde{X}^c(t_m)) \, f_{m+1}(\tilde{X}^c(t_{m+1}))] \\ &= \mathbb{E}[f_1(\tilde{X}^c(t_1)) \cdots f_m(\tilde{X}^c(t_m)) \, \Psi_{t_{m+1} - t_m} f_{m+1}(\tilde{X}^c(t_m))], \end{aligned}$$

and the Markov property is established. $\square$

REMARK 3.14. $\Psi$ defined in (3.22) is the semigroup of $\tilde{X}^c$. Lemma 3.13 implies that $\Psi_t f(\cdot)$ is continuous for any $t \geq 0$ and $f : \bar{\mathbb{N}}^* \to \mathbb{R}$ continuous. We thus have that $\tilde{X}^c$ is a Feller process. This and the Markov property just proved imply the strong Markov property of $\tilde{X}^c$.

Next follows a result establishing in particular the lack of strong continuity of the K-process with $c = 0$.

LEMMA 3.15. For every $y \in \bar{\mathbb{N}}^*$, we have that $\mathbb{P}(\tilde{X}^{0,y}(t) = \infty) = 0$ for every $t > 0$.

REMARK 3.16. The statement of Lemma 3.15 does not hold for $c > 0$. In this case, it can actually be shown that $\mathbb{P}(\tilde{X}^c(t) = \infty) > 0$ for every $t > 0$, $y \in \bar{\mathbb{N}}^*$. It can also be shown that the process is strongly continuous in this case.

REMARK 3.17. We note that $\tilde{X}^{c,\infty}(0) = \infty$ almost surely for every $c \geq 0$.

PROOF OF LEMMA 3.15. For $m \geq 1$, and $t > 0$, let $\theta_{m,t}$ be the time spent by $\tilde{X}^{0,y}$ outside $\{1, \ldots, m\}$ up to time $t$. Clearly

$$(3.29) \qquad \theta_{m,t} = \sum_{x > m} \gamma(x) \sum_{i=1}^{N_{\Xi_t}^{(x)}} T_i^{(x)},$$

where $\Xi$ is the inverse function of $\Gamma$. It is also clear that

$$(3.30) \qquad \int_0^t 1_{\{\infty\}}(\tilde{X}^{0,y}(s)) \, ds \leq \theta_{m,t}$$

for every $m \geq 1$ and $t > 0$, where $1.$ is the usual indicator function, and that

$$\theta_{m,t} \to 0 \tag{3.31}$$

almost surely as $m \to \infty$ for every $t > 0$. Thus the left-hand side of (3.30) vanishes almost surely and dominated convergence implies that

$$\mathbb{E}\bigg(\int_0^t 1_{\{\infty\}}(\tilde{X}^{0,y}(s))\,ds\bigg) = \int_0^t \mathbb{P}(\tilde{X}^{0,y}(s) = \infty)\,ds = 0 \tag{3.32}$$

for every $t$. This proves the assertion of the lemma for Lebesgue-almost every $t$. The Markov property of $\tilde{X}^{0,y}$ can now be used to extend the result to every $t$. □

We close this section with a computation related to the Green function of $\tilde{X}^c$; this will lead to an identification of $\tilde{X}^c$ above and $X^c$ defined in Section 2.3.

Let $\tau^{\{x\}} = \inf\{t\,;\,\tilde{X}^c(t) = x\}$. We have that under $\mathbb{P}_\infty$

$$\tau^{\{x\}} = \Gamma_c^{(x)}(\sigma_1^{(x)}), \tag{3.33}$$

where for $x \in \mathbb{N}^*$, $s \geq 0$

$$\Gamma_c^{(x)}(s) = \sum_{y \neq x} \gamma_y \sum_{i=1}^{N_s^{(y)}} T_i^{(y)} + cs. \tag{3.34}$$

It is now straightforward to compute the Laplace transform of $\tau^{\{x\}}$ for the process started at $\infty$. We obtain

$$\mathbb{E}_\infty(e^{-\lambda \tau^{\{x\}}}) = \mathbb{E}_\infty(e^{-\lambda \Gamma_c^{(x)}(\sigma_1^{(x)})}) = \int_0^\infty \mathbb{E}(e^{-\lambda \Gamma_c^{(x)}(s)})e^{-s}\,ds$$

$$= \int_0^\infty \mathbb{E}\bigg(\exp\bigg\{-\lambda \sum_{y \neq x} \gamma_y \sum_{i=1}^{N_s^{(y)}} T_i^{(y)}\bigg\}\bigg)e^{-(1+c)s}\,ds, \tag{3.35}$$

for $\lambda > 0$, where in the second equality we have used the independence of $\sigma_1^{(x)}$ and (the random variables in) $\sum_{y \neq x} \gamma_y \sum_{i=1}^{N_s^{(y)}} T_i^{(y)}$. We leave it as an exercise to compute the expectation inside the integral in (3.35), and to conclude that the integral equals

$$\bigg(1 + c\lambda + \sum_{y \neq x} \frac{\lambda \gamma_y}{1 + \lambda \gamma_y}\bigg)^{-1}. \tag{3.36}$$

We note that this expression is the same as that for the corresponding transform for $X^c$ in Section 2.

Now, since the only transitions are from states in $\mathbb{N}^*$ to $\infty$ and back, we have a decomposition as in (2.12) for the Green kernel of $\tilde{X}^c$ starting at $\infty$, and we get the case of general initial condition from the case of $\infty$ initial condition as in the computation in Section 2. We readily conclude from the remark at the end of the previous paragraph that the Green functions of $\tilde{X}^c$ and $X^c$ coincide, and thus, since these are both càdlàg Markov processes, they must have the same distribution for any initial law.

**4. A characterization result.** The striking property of K-processes that at the entrance time of the process in finite subsets (starting from outside) the distribution is uniform (see Remark 2.9) leads to a natural question: which other processes have this property? Below we see that, under natural assumptions, the answer is *none*; that is, that property characterizes K-processes.

THEOREM 4.1. *Let $\gamma$ be as in (2.2) and $Y = (Y(t), t \geq 0)$ be a process on $\bar{\mathbb{N}}^*$ with the following four properties:*

  (i) *$Y$ is càdlàg.*
  (ii) *$Y$ is strong Markov.*
  (iii) *Starting from any point $i \in \mathbb{N}^*$, $Y$ waits for an exponential time of mean $\gamma(i)$ before jumping.*
  (iv) *Starting from $\infty$, for any finite $A \subset \mathbb{N}^*$, we have $\tau_A < \infty$ almost surely, where*

$$(4.1) \qquad \tau_A = \inf\{t \geq 0 : Y(t) \in A\},$$

*with $\inf \varnothing = \infty$, and the law of $Y(\tau_A)$ is uniform on $A$.*

*Then, $Y$ is a K-process with parameters $\gamma$ and $c$, for some $c \geq 0$.*

REMARK 4.2. Fukushima and collaborators, as well as other authors, have recently studied one-point extensions of certain Markov processes beyond a killing time (see [11, 18] and references therein). The K-processes can be viewed as one-point extensions of processes in $\mathbb{N}^*$ that are killed after the first jump. With this point of view, and even though the K-processes do not satisfy some of the conditions in the above references (like Condition (A.2) in [11]; another condition would require $c = 0$ in our case), Theorem 4.1 is similar (in its particular context) to their results. But there is an important difference in that, while they depart from a reversibility condition (more generally, a duality condition) with respect to an excessive measure for the process, we have a condition on the jump rates and entrance laws. It is nevertheless remarkable that entrance laws play a crucial role in their approach (it also could be said that for the K-processes the jump rates are directly related to a stationary measure for the process).

PROOF OF THEOREM 4.1. The strategy is to consider the process restricted to $\{1, \ldots, n, \infty\}$, and show that it must have the same distribution as $\tilde{X}_n^c$. The result then follows by taking $n \to \infty$.

We start by showing that from any state in $\mathbb{N}^*$, $Y$ jumps to $\infty$ almost surely. Let $i, j \in \mathbb{N}^*$ be such that $i \neq j$, and for $n \geq i \vee j$ let

$$(4.2) \qquad A_n = \{1, \ldots, n\}.$$

Then

$$(4.3) \quad \{Y(\tau_{A_{n+1}}) = j\} \cup \{Y(\tau_{A_{n+1}}) = i, Y(\tau'_{n+1}) = j\} \subset \{Y(\tau_{A_{n+1} \setminus \{i\}}) = j\},$$

where

$$(4.4) \qquad \tau'_n = \inf\{t \geq \tau_{A_n} : Y(t) \neq Y(\tau_{A_{n+1}})\}.$$

Thus,

$$(4.5) \qquad \begin{aligned} &\mathbb{P}_\infty(Y(\tau_{A_{n+1}}) = j) + \mathbb{P}_\infty(Y(\tau_{A_{n+1}}) = i, Y(\tau'_{n+1}) = j) \\ &\qquad \leq \mathbb{P}_\infty(Y(\tau_{A_{n+1} \setminus \{i\}}) = j), \end{aligned}$$

and using (ii)–(iv)

$$(4.6) \qquad \frac{1}{n+1} + \frac{1}{n+1} p_{ij} \leq \frac{1}{n},$$

where $p_{ij}$ is the transition probability from $i$ to $j$. It follows that $p_{ij} \leq 1/n$, and since $n$ can be taken arbitrarily large, we conclude that $p_{ij} = 0$, and the claim at the beginning of the paragraph follows.

Now let us consider the process obtained from $Y$ by suppressing jumps outside $\bar{A}_n := A_n \cup \{\infty\}$ [see (4.2)]. Let us call it $Y_n$. More precisely, let

$$\mathcal{A}_n(t) = \int_0^t \mathbf{1}_{\bar{A}_n}(Y(s)) \, ds \quad \text{and} \quad Y_n(t) = Y(\mathcal{A}_n^{-1}(t)),$$

where $\mathbf{1}_A$ is the usual indicator function of a set $A$, and $\mathcal{A}_n^{-1}$ is the right-continuous inverse of $\mathcal{A}_n$.

It is readily seen that $\mathcal{A}_n(t) \uparrow t$ as $n \uparrow \infty$ uniformly in $t \leq T$ for every $T$.

We can also prove that

$$(4.7) \qquad Y_n \to Y \qquad \text{as } n \to \infty \text{ in Skorohod space.}$$

This can be argued through parallel steps to those in the proof of Lemma 3.11, with $\mathcal{A}_n^{-1}$ replacing $\Gamma_n$, and the $S_i^m$'s replaced by the successive exit times of $\{1, \ldots, m\}$ by $Y_n$. Notice that (iv) implies the existence of $n_m$ as in that proof.

We observe that $Y_n$ also satisfies (i) and (ii) (see [23], Theorem 65.9). Starting at $i \leq n$, it waits an exponential time of mean $\gamma(i)$ and then jumps.

The state space of $Y_n$ may be either $\bar{A}_n$ or $A_n$. [The latter possibility happens if $Y$ is the K$(\gamma, 0)$-process; note that $\int_0^\infty \mathbf{1}_{\{\infty\}}(Y(s))\,ds = 0$ almost surely in that case—see the proof of Lemma 3.15 for an argument.] If the latter case happens, then $Y_n$ is a continuous-time Markov chain on $A_n$ satisfying (i)–(iii). To completely characterize it, we need only determine the transition probabilities. But property (iv) of $Y$ implies that these must be uniform, that is, $p_{ij} \equiv 1/n$. This means that $Y_n$ is equidistributed with $\bar{X}_n^0$ defined in (3.5). Now this, Lemma 3.11 and (4.7) imply that $Y$ is equidistributed with $\bar{X}^0$.

It remains to consider the case where the state space of $Y_n$ is $\bar{A}_n$. In this case $Y_n$ clearly also satisfies (iii)–(iv). We need only determine the mean holding time at $\infty$, say $\gamma_n(\infty)$. For that we reason as follows.

We can obtain $Y_n$ by suppressing jumps of $Y_{n+1}$ outside $\{1, \ldots, n\} \cup \{\infty\}$. Since upon leaving $\infty$, the process $Y_{n+1}$ has probability $1/(n+1)$ to jump to $n + 1$, we see that the holding time at $\infty$ in $Y_n$ can be seen as a sum of independent holding times at $\infty$ in $Y_{n+1}$. The number of terms in the sum is a geometric random variable with success parameter $1/(n+1)$, independent of the holding times at $\infty$ in $Y_{n+1}$. We conclude that

$$\gamma_n(\infty) = \gamma_{n+1}(\infty)\frac{n+1}{n},$$

and thus that $\gamma_n(\infty) = c/n$ for some constant $c \geq 0$. We then see that $Y_n$ is equidistributed with $\bar{X}_n^c$ defined in (3.5) for every $n \geq 1$, and the conclusion that $Y$ is equidistributed with $\bar{X}^c$ follows exactly as above.    $\square$

## 5. A scaling limit for the trap model in the complete graph.

The trap model in the complete graph [7] can be described as a continuous-time symmetric Markov chain $Y_n = \{Y_n(t), t \geq 0\}$ in the complete graph $K_n$ with $n$ vertices such that the average holding times $\tau := \{\tau_x, x \in K_n\}$ are an i.i.d. family of positive random variables equidistributed with a r.v. $\tau_0$ which is in the basin of attraction of a stable law of degree $\alpha < 1$, that is,

$$(5.1) \qquad \mathbb{P}(\tau_0 > t) = \frac{L(t)}{t^\alpha}, \qquad t > 0,$$

where $L$ is a slowly varying function at infinity.

We will show in this section that in an appropriate sense, in an appropriate time scale, $Y_n$ converges in distribution as $n \to \infty$ to a K-process with $c = 0$.

We start by identifying the vertices of $K_n$ with $A_n = \{1, \ldots, n\}$ for all $n \geq 1$, in such a way that $\{\tau_i^{(n)}, 1 \leq i \leq n\}$ is in decreasing order, that is, $(\tau_1^{(n)}, \ldots, \tau_n^{(n)})$ is the reverse order statistics of $(\tau_1, \ldots, \tau_n)$, an i.i.d. sample of size $n$ of $\tau_0$.

We can describe $Y_n$ then as a continuous-time Markov chain in $A_n$ with mean holding time at $i \in A_n$ given by $\tau_i^{(n)}$ and uniform in $A_n$ transition probabilities for all starting point $i \in A_n$.

Let us view $(\tau_1^{(n)}, \ldots, \tau_n^{(n)})$ as a random measure $\gamma_n$ on $\mathbb{N}^*$ such that

$$(5.2) \qquad \gamma_n(\{i\}) = \begin{cases} \tau_i^{(n)}, & \text{if } i \in A_n, \\ 0, & \text{otherwise.} \end{cases}$$

We refer the reader to Section 3 of [16] for more on the context of the next result. We present the main points below.

Consider the increasing Lévy process $V_x$, $x \in \mathbb{R}$, $V_0 = 0$, with stationary and independent increments given by

$$(5.3) \qquad \mathbb{E}[e^{ir(V_{x+x_0} - V_{x_0})}] = e^{\alpha x \int_0^\infty (e^{irw} - 1) \, w^{-1-\alpha} \, dw}$$

for any $x_0 \in \mathbb{R}$ and $x \geq 0$. Let $\rho$ be the (random) Lebesgue–Stieltjes measure on the Borel sets of $\mathbb{R}$ associated to $V$, that is,

$$(5.4) \qquad \rho((a, b]) = V_b - V_a, \qquad a, b \in \mathbb{R}, \ a < b.$$

Then

$$(5.5) \qquad d\rho = dV = \sum_j w_j \, \delta(x - x_j),$$

where the (countable) sum is over the indices of an inhomogeneous Poisson point process $\{(x_j, w_j)\}$ on $\mathbb{R} \times (0, \infty)$ with density $dx \, \alpha w^{-1-\alpha} \, dw$.

Let $\gamma = \{\gamma(i), i \in \mathbb{N}^*\}$ denote the weights of $\rho$ in $[0, 1]$ in decreasing order, that is, making $\mathcal{R} = \{\gamma(\{x\}), x \in [0, 1]\}$,

$$(5.6) \quad \gamma(1) = \max \mathcal{R}; \qquad \gamma(i) = \max[\mathcal{R} \setminus \{\gamma(1), \ldots, \gamma(i-1)\}], \qquad i \geq 2.$$

REMARK 5.1. $\gamma$ thus defined almost surely satisfies the conditions on the paragraph of (2.2).

Let now

$$(5.7) \qquad c_n = (\inf\{t \geq 0 : \mathbb{P}(\tau_0 > t) \leq n^{-1}\})^{-1}.$$

and $\tilde{\gamma}_n = c_n \gamma_n$, that is, $\tilde{\gamma}_n$ is a (random) measure in $\mathbb{N}^*$ such that

$$(5.8) \qquad \tilde{\gamma}_n(\{i\}) = \begin{cases} c_n \tau_i^{(n)}, & \text{if } i \in A_n, \\ 0, & \text{otherwise.} \end{cases}$$

THEOREM 5.2. Let $\tilde{Y}_n$ be the process in $A_n$ such that for $t \geq 0$, $\tilde{Y}_n(t) \equiv Y_n(c_n^{-1}t)$. Suppose $\tilde{Y}_n(0) \equiv \mathcal{Y}_n$ converges weakly to a random variable $\mathcal{Y}$ in $\bar{\mathbb{N}}^*$. Then, as $n \to \infty$,

$$(5.9) \qquad (\tilde{Y}_n, \tilde{\gamma}_n) \Rightarrow (Y, \gamma),$$

where, given $\gamma$, $Y$ is a $K(\gamma, 0)$-process with $Y(0)$ distributed as $\mathcal{Y}$, and $\Rightarrow$ denotes weak convergence in the product of the Skorohod topology and the vague topology in the space of finite measures on $\bar{\mathbb{N}}^*$.

REMARK 5.3. Given our choice of $\bar{\mathbb{N}}^*$ as state space of both $Y_n$ and $Y$ in the above convergence result, the ordering of $\tau$ and $\gamma$ is a natural imposition: in this way we have that the $\tau_i$'s that matter in the chosen scaling limit get naturally assigned to the points of the state space. (Without the ordering of $\tau$, the locations of the $\tau_i$'s that do not vanish in the limit go to infinity.)

Other choices of state space would not require the ordering. For example, one could rescale the space variable by $n$, and take the state space as the interval $[0, 1]$; then the pairs (rescaled $\tau$'s, their rescaled locations) would converge in distribution to the Poisson point process alluded to above [right below (5.5)] in $[0, 1] \times (0, \infty)$.

REMARK 5.4. As is apparent from the proof of Theorem 5.2, the result remains valid whenever the rescaled (and reordered) $\tau$ (its distribution could also depend on $n$), possibly under a different scaling, converges to some $\gamma$ satisfying the conditions on the paragraph of (2.2). In the case of a different scaling for $\tau$, the time of $\tilde{Y}$ should of course be rescaled likewise.

PROOF OF THEOREM 5.2. We may assume that $\mathcal{Y}_n \to \mathcal{Y}$ as $n \to \infty$ almost surely. Following the strategy in Section 3 of [16], we will couple $(\tilde{Y}_n, \tilde{\gamma}_n)$ to $(Y, \gamma)$ and establish (5.9) as a strong convergence.

For $i \in A_n$, let

$$(5.10) \qquad \bar{\tau}_i^{(n)} = \frac{1}{c_n} g_n(V_{i/n} - V_{(i-1)/n}),$$

where $g_n$ is defined as follows. Let $G : [0, \infty) \to [0, \infty)$ satisfy

$$(5.11) \qquad \mathbb{P}(V_1 > G(x)) = \mathbb{P}(\tau_0 > x) \qquad \text{for all } x \geq 0$$

and let $g_n : [0, \infty) \to [0, \infty)$ be defined as

$$(5.12) \qquad g_n(x) = c_n G^{-1}(n^{1/\alpha} x) \qquad \text{for all } x \geq 0.$$

Let $\{\hat{\tau}_i^{(n)}, i \in A_n\}$ be $\{\bar{\tau}_i^{(n)}, i \in A_n\}$ in decreasing order, and

$$(5.13) \qquad \hat{\gamma}_n(\{i\}) = \begin{cases} c_n \hat{\tau}_i^{(n)}, & \text{if } i \in A_n, \\ 0, & \text{otherwise.} \end{cases}$$

It readily follows from Proposition 3.1 in [16] that $\hat{\gamma}_n$ and $\tilde{\gamma}_n$ have the same distribution for every $n$, and that almost surely

$$(5.14) \qquad \hat{\gamma}_n \to \gamma \qquad \text{as } n \to \infty$$

[where the first $\to$ in (5.14) means vague convergence].

Let now $\mathcal{N}$ and $\mathcal{T}$ be as in Section 3. For $t \geq 0$, let

$$\hat{\Gamma}_n(t) = \hat{\gamma}_n(\mathcal{Y}_n) \, T_0 + \sum_{x=1}^{n} \hat{\gamma}_n(x) \sum_{i=1}^{N_t^{(x)}} T_i^{(x)}, \tag{5.15}$$

$$\hat{\Gamma}(t) = \gamma(\mathcal{Y}) \, T_0 + \sum_{x=1}^{\infty} \gamma(x) \sum_{i=1}^{N_t^{(x)}} T_i^{(x)}, \tag{5.16}$$

where we write $\hat{\gamma}_n(x)$ and $\gamma(x)$ for $\hat{\gamma}_n(\{x\})$ and $\gamma(\{x\})$, respectively, and

$$\hat{Y}_n(t) = \begin{cases} \mathcal{Y}_n, & \text{if } 0 \leq t < \hat{\gamma}_n(\mathcal{Y}_n) \, T_0, \\ x, & \text{if } \hat{\Gamma}_n(\sigma_j^{(x)}-) \leq t < \hat{\Gamma}_n(\sigma_j^{(x)}) \\ & \quad \text{for some } 1 \leq x \leq n, \, j \geq 1. \end{cases} \tag{5.17}$$

$$\hat{Y}(t) = \begin{cases} \mathcal{Y}, & \text{if } 0 \leq t < \gamma(\mathcal{Y}) \, T_0, \\ x, & \text{if } \hat{\Gamma}(\sigma_j^{(x)}-) \leq t < \hat{\Gamma}(\sigma_j^{(x)}) \\ & \quad \text{for some } 1 \leq x < \infty, \, j \geq 1 \\ \infty, & \text{otherwise.} \end{cases} \tag{5.18}$$

See (3.1), (3.2), (3.4), (3.5) above. One readily checks that $(\hat{Y}_n, \hat{\gamma}_n)$ has the same distribution as $(Y_n, \hat{\gamma}_n)$ for every $n \geq 1$ (see Proposition 3.1 of [16]).

We claim now that

$$\hat{Y}_n \to Y \qquad \text{as } n \to \infty \tag{5.19}$$

almost surely in Skorohod space.

The proof of (5.19) is similar to that of Lemma 3.11, with modifications to account for a dependence of $\hat{\gamma}_n$ on $n$. Equation (5.14) is of course crucial. We indicate the main steps next.

For $n \geq y$, $m \in \mathbb{N}^*$, let $\delta_m$, $\{S_1^m < S_2^m < \cdots\}$ and $L_n^m$ be as in that proof. We now have that $\min_{0 \leq i \leq L_n^m - 1} [\hat{\Gamma}_n(S_{i+1}^m-) - \hat{\Gamma}_n(S_i^m)] > 0$ almost surely for $n \geq 1$. Define next $\hat{\lambda}_n^m : [0, \hat{\Gamma}_n(S_{L_n^m}^m)] \to \mathbb{R}^+$ as follows:

$$\hat{\lambda}_n^m(t) = \frac{\gamma(\mathcal{Y})}{\hat{\gamma}_n(\mathcal{Y}_n)} t \qquad \text{if } 0 \leq t < \hat{\gamma}_n(\mathcal{Y}_n) \, T_0, \tag{5.20}$$

and for $0 \leq i \leq L_n^m - 1$ and $\hat{\Gamma}_n(S_i^m) \leq t \leq \hat{\Gamma}_n(S_{i+1}^m)$, let

$$\hat{\lambda}_n^m(t) = \begin{cases} \hat{\Gamma}(S_i^m) + \dfrac{\hat{\Gamma}(S_{i+1}^m-) - \hat{\Gamma}(S_i^m)}{\hat{\Gamma}_n(S_{i+1}^m-) - \hat{\Gamma}_n(S_i^m)} [t - \hat{\Gamma}_n(S_i^m)], \\ \qquad\qquad\qquad\qquad \text{if } \hat{\Gamma}_n(S_i^m) \leq t \leq \hat{\Gamma}_n(S_{i+1}^m-), \\ \hat{\Gamma}(S_{i+1}^m-) + \dfrac{\hat{\Gamma}(S_{i+1}^m) - \hat{\Gamma}(S_{i+1}^m-)}{\hat{\Gamma}_n(S_{i+1}^m) - \hat{\Gamma}_n(S_{i+1}^m-)} [t - \hat{\Gamma}_n(S_{i+1}^m-)], \\ \qquad\qquad\qquad\qquad \text{if } \hat{\Gamma}_n(S_{i+1}^m-) \leq t \leq \hat{\Gamma}_n(S_{i+1}^m). \end{cases} \tag{5.21}$$

It has the following properties. For all $T > 0$, $m \in \mathbb{N}^*$ and $n \geq m \vee y$

$$
\sup_{0 \leq t \leq T} |\hat{\lambda}_n^m(t) - t|
$$

(5.22)

$$
\leq \max\{|\hat{\Gamma}(S_i^m-) - \hat{\Gamma}_n(S_i^m)|, |\hat{\Gamma}(S_i^m) - \hat{\Gamma}_n(S_i^m)|; 0 \leq i \leq L_n^m\},
$$

where $\hat{\Gamma}_n(0-) \equiv \hat{\Gamma}(0-) = 0$, and the right-hand side of (5.22) vanishes almost surely as $n \to \infty$. Furthermore,

(5.23)

$$
\sup_{0 \leq t \leq T} \mathrm{dist}(\hat{Y}(\hat{\lambda}_n^m(t)), \hat{Y}_n(t)) \leq \delta_m,
$$

since for $t \in [0, T]$, $\hat{Y}(\hat{\lambda}_n^m(t))$ and $\hat{Y}_n(t)$ coincide when either one is in $\{1, \ldots, m\}$.

The remainder of the argument follows along the exact same lines as those in the proof of Lemma 3.11.  $\square$

REMARK 5.5.  It is natural to ask about the marginal distribution of the limiting process $Y$ of Theorem 5.2. Even though we do not study this process at length here, one thing that can be readily checked is that it is not Markovian, since in particular the distribution of the waiting times is not exponential (but a mixture thereof), and the waiting times at different sites are also not independent (due to the dependence of the distribution on the mixture kernels of different sites).

One positive thing that can be said of $Y$ is that it exhibits aging. See Theorem 5.11 and Remark 5.12, which obviously holds for Corollary 5.15 as well.

5.1. *Aging.*  Aging results can be viewed as scaling limits for averaged two-time correlation functions of a given dynamics of a disordered system. The averaging is with respect to the disorder distribution. The system should be started at high temperature, and then abruptly cooled down, evolving thence on at low temperature. Loosely speaking, aging would amount to the following. Given a dynamics described by the process $X$ with a disordered set of parameters $\tau$, the following would be an aging result:

(5.24)

$$
\lim_{\substack{t, t' \to \infty \\ t'/t \to \theta}} \mathbb{E}\{\mathbb{E}_\mu[\Phi(t, t'; X)|\tau]\} = \mathcal{R}(\theta),
$$

where $\mu$ is a measure on state space; $\mathbb{E}_\mu(\cdot|\tau)$ indicates the expectation with respect to $X$ with initial distribution given by $\mu$, with parameters fixed at $\tau$; $\Phi(t, t'; X)$ is a function of the piece of trajectory $X([t, t + t']) = \{X(s), s \in [t, t + t']\}$; and $\mathcal{R}$ is a nontrivial function of real scaling factor $\theta > 0$. The initial distribution $\mu$ should reflect a high temperature, and the distribution of the parameters, a low temperature. See [8] and references therein.

For a mean-field model like the trap model in the complete graph, there is a volume dependence, and one must take the infinite volume limit $(n \to \infty)$; that should be done before or together with the time limit. The former is done in [7] for

$$(5.25) \qquad \Phi_1(t, t'; X) = 1\{X(s) = X(t), \, s \in [t, t + t']\}.$$

$\mu = \mu_n$ is taken uniformly distributed in $\{1, \dots, n\}$, reflecting the high temperature of the initial state, and the tail parameter $\alpha < 1$ corresponds to the low temperature thence prevailing.

One other function that is often considered is

$$(5.26) \qquad \Phi_2(t, t'; X) = 1\{X(t) = X(t + t')\}.$$

One could also take the volume and time limits together, using the scaling limit of Theorem 5.2. For that let us suppose that, for all $t, t' > 0$, $\Phi(t, t'; \cdot)$ is almost surely continuous [with respect to the distribution of $(Y, \gamma)$]. Then, by Theorem 5.2,

$$(5.27) \qquad \lim_{n \to \infty} \mathbb{E}\{\mathbb{E}_{\mu_n}[\Phi(t, t'; \tilde{Y}_n)|\tau]\} = \mathbb{E}\{\mathbb{E}_\infty[\Phi(t, t'; Y)|\gamma]\},$$

with $Y(0) = \infty$, where for $x \in \bar{\mathbb{N}}^*$, $\mathbb{E}_x(\cdot|\gamma)$ denotes the expectation with respect to the distribution of $Y$ started at $x$, with parameters fixed at $\gamma$. $t, t' > 0$ are now macroscopic. On those times the dynamics is already close enough to equilibrium to disallow aging: the right-hand side of (5.27) is *not* a function of the ratio $t'/t$ only. To find aging, we should move away from equilibrium, by taking the further limit as $t, t' \to 0$ while $t'/t \to \theta > 0$. We then say that aging takes place in this context (for both the trap model in the complete graph and the limiting disordered K-process) if

$$(5.28) \qquad \lim_{\substack{t, t' \to 0 \\ t'/t \to \theta}} \mathbb{E}\{\mathbb{E}_\infty[\Phi(t, t'; Y)|\gamma]\} = \mathcal{R}'(\theta)$$

exists and $\mathcal{R}'$ is nontrivial.

REMARK 5.6. In taking the volume and time limits as in [7], one enters what could be termed a (long) microscopic-time aging regime for the trap model in the complete graph, while the latter way of taking those limits gets one in a (short) macroscopic-time aging regime. Our next result indicates that, at least as far as $\Phi_1$ is concerned, the two regimes agree.

REMARK 5.7. Instead of scaling time as in (5.27), namely with the scale of the largest $\tau_x$'s in $K_n$, in view of the further limit $\lim_{t, t' \to 0; \, t'/t \to \theta}$, it is natural to use a lower-order (divergent) scaling. This could be termed a *mesoscopic* aging regime, and it is the approach of [3] to establishing aging for the trap model in the complete graph. As far as $\Phi_1$ is concerned, the mesoscopic aging regime agrees with the microscopic and macroscopic regimes; see [3].

Next we state an aging result for $\Phi$ in a certain class of functions including the usual examples $\Phi_1$ and $\Phi_2$ and satisfying some continuity and spatial homogeneity conditions (which seem natural if one sees this as a mean-field model), with no intention at full generality, however. Let $\Pi$ be the space of càdlàg paths on $\bar{\mathbb{N}}^*$, and consider the class of functions $\Phi : \mathbb{R}^+ \times \mathbb{R}^+ \times \Pi \to \mathbb{R}$ with the following properties:

$$\Phi(t, t'; \zeta) = \Phi(t, t'; \zeta([t, t + t'])), \tag{5.29}$$

where $\zeta([t, t + t'])$ is $\zeta$ restricted to $[t, t + t']$, with the scaling property: for all $t, t' > 0$,

$$\Phi(t, t'; \zeta) = \Phi(1, t'/t; \zeta^t), \tag{5.30}$$

where $\zeta^t(\cdot) = \zeta(t \cdot)$. Notice that $\Phi_1$ and $\Phi_2$ above have this property. Consider now the following path segments: for $\theta > 0$, $x \in \bar{\mathbb{N}}^*$: $\eta_{x,\theta} = \eta_{x,\theta}([1, 1 + \theta)) \equiv x$; $\bar{\eta}_{x,\theta} = \eta_{x,\theta}([1, 1 + \theta]) \equiv x$. We make the following further assumptions on $\Phi$:

$$\Phi(1, \theta; \bar{\eta}_{x,\theta}) = \Psi_1(\theta) \qquad \forall x \in \mathbb{N}^*, \tag{5.31}$$

for some real function $\Psi_1$, that is, $\Phi(1, \theta; \bar{\eta}_{x,\theta})$ does not depend on $x$ for finite $x$; for $0 < s < \theta$,

$$\Phi(1, \theta; \eta_{x,s} \circ \eta) = \Xi(x, s, \eta) \qquad \forall x \in \mathbb{N}^* \tag{5.32}$$

for all segment $\eta = \eta([1 + s, 1 + \theta])$ in $[1 + s, 1 + \theta]$ of a path in Skorohod space with $\eta(1 + s) \neq x$, where $\circ$ stands for concatenation, and $\Xi$ is a given function with the following properties. We first give some definitions. For $r > 0$, let $\Pi_u$ be the space of càdlàg paths in $\bar{\mathbb{N}}^*$ of length $u$, and for $v > 0$ fixed, let $\mathcal{X}_v = \bigcup_{u \in [0,v]} (\{u\} \times \Pi_{v-u})$. Let now $\theta > 0$ be fixed. We then have that $\Xi : \bar{\mathbb{N}}^* \times \mathcal{X}_\theta \to \mathbb{R}$ is such that:

(i)
$$\Xi \text{ is uniformly bounded;} \tag{5.33}$$

(ii)
$$\text{for all } \eta \in \Pi_{\theta-s}, \Xi(x, s, \eta) = \Xi(y, s, \eta) \text{ whenever } x, y \notin \eta; \tag{5.34}$$

(iii)
for all $x \in \mathbb{N}^*$, the function $\eta \to \Xi(x, s, \eta)$ is continuous *in the sup*

$$\textit{norm on } \Pi_{\theta-s} \textit{ for } \eta = \hat{\eta}_{\infty,\theta-s} \in \Pi_{\theta-s} \textit{ with } \hat{\eta}_{\infty,\theta-s} \equiv \infty, \tag{5.35}$$

uniformly in $0 < s < \theta$.

REMARK 5.8. $\Phi_1$ and $\Phi_2$ given in (5.25), (5.26) satisfy (5.29)–(5.35). Other examples can be obtained by taking $\Xi(x, s, \eta) = f(s)$ for all $x \in \mathbb{N}^*$, $\eta \in \Pi_{\theta-s}$, where $f$ is any continuous function in $[0, \theta]$.

REMARK 5.9. The uniformity assumption in (5.35) above is for simplicity. See Remarks 5.18 and 5.20 below.

REMARK 5.10. The lack of dependence on finite $x$ assumed in both (5.31) and (5.34) is not artificial if one takes into account that the model where $\Phi$ will be measured is mean-field, and thus the space coordinate is not relevant. The distinction between finite $x$ and infinity is nevertheless desirable.

For $x \in \mathbb{N}^*$ and $0 < s < \theta$, let $\Psi_2(s,\theta) = \Xi(x, s, \hat{\eta}_{\infty,\theta-s})$. Notice that the latter function does not depend on $x \in \mathbb{N}^*$ by assumption (ii). We make the following assumptions on $\Psi_2$, for simplicity: for all $\theta > 0$

$$(5.36) \quad \Psi_2(0,\theta) := \lim_{s \downarrow 0} \Psi_2(s,\theta), \qquad \Psi_2(\theta,\theta) := \lim_{s \uparrow \theta} \Psi_2(s,\theta) \qquad \text{exist}$$

and

$$(5.37) \qquad \Psi_2'(s,\theta) := \frac{d}{ds}\Psi_2(s,\theta) \in L_1([0,\theta],\,dx).$$

We can now state the main results of this subsection. For simplicity, we make $t' = \theta t$. We start with a particular case.

THEOREM 5.11. *For $\gamma$ as in* (5.3)–(5.6), *$\Phi_1$ as in* (5.25), *and $t, \theta > 0$, let*

$$(5.38) \qquad \Lambda_t(\theta) = \mathbb{E}_\infty[\Phi_1(t, \theta t; Y)|\gamma].$$

*Then almost surely for every $\theta > 0$*

$$(5.39) \qquad \lim_{t \to 0} \Lambda_t(\theta) = \Lambda(\theta),$$

*where $\Lambda$ is a (nontrivial) function to be exhibited below [see* (5.95), *and also Proposition* 5.22*].*

REMARK 5.12. This is an almost sure aging result. The averaged form of (5.27) follows by dominated convergence.

REMARK 5.13. As anticipated in Remark 5.6 above, $\Lambda$ is the same as the one obtained in [7] by taking limits in a different order and in a different way (see the discussion before and up to that remark).

REMARK 5.14. Were $Y$ a $K(\gamma, c)$-process with $c > 0$, then the limit in (5.39) would be a trivial one. Indeed, $Y$ would be strongly continuous in this case, meaning that $\mathbb{P}_\infty(Y_t = \infty) \to 1$ as $t \to 0$ (see the definition of "strongly continuous" in Remark 3.2, and see Remark 3.16 above). We would thus have $\Lambda \equiv 1$. See (end of) Remark 3.2 above.

COROLLARY 5.15. *Let $\Phi$ be as in* (5.29)–(5.37). *If* (5.39) *holds, then*

$$\begin{aligned}
(5.40) \qquad & \lim_{t \to 0} \mathbb{E}_\infty[\Phi(t, \theta t; Y) | \gamma] \\
& = \Psi_2(0, \theta) + [\Psi_1(\theta) - \Psi_2(\theta, \theta)] \Lambda(\theta) + \int_0^\theta \Psi_2'(s, \theta) \Lambda(s) \, ds.
\end{aligned}$$

REMARK 5.16. For the above result, we may take $\gamma$ fixed such that (5.39) holds [as well as the assumptions on the paragraph of (5.6)].

REMARK 5.17. For both $\Phi_1$ and $\Phi_2$ [see (5.25), (5.26)], we have $\Psi_1 \equiv 1$ and $\Psi_2 \equiv 0$, so, from (5.40), $\Lambda(\theta)$ is their common aging limit.

REMARK 5.18. $\Lambda$ turns out to be continuously differentiable in $[0, \infty)$; we can thus integrate by parts in the right-hand side of (5.40) to obtain that

$$(5.41) \qquad \lim_{t \to 0} \mathbb{E}_\infty[\Phi(t, \theta t; Y) | \gamma] = \Psi_1(\theta) \Lambda(\theta) - \int_0^\theta \Psi_2(s, \theta) \Lambda'(s) \, ds,$$

where $\Lambda'(s) = \frac{d}{ds} \Lambda(s)$. For the result in this form we do not require the uniformity assumption in (5.35), nor assumptions (5.36), (5.37). See Remark 5.20 below.

REMARK 5.19. In the proof of Corollary 5.15 below, we will use the fact that for each $\gamma$ satisfying the conditions of the paragraph of (2.2)—in particular, for each $\gamma$ in a full measure event; see Remark 5.1 above—and all $t > 0$, the distribution of $Y^t$ given $\gamma$ is the same as that of $Y$ given $\gamma^t := t^{-1}\gamma$. This follows immediately from the definition of K-processes (see Definition 3.1 and preceding paragraphs). We thus have that for all such $\gamma$, and all bounded measurable function $F$ on Skorohod space,

$$(5.42) \qquad \mathbb{E}[F(Y^t) | \gamma] = \mathbb{E}[F(Y) | \gamma^t].$$

PROOF OF COROLLARY 5.15. Consider the conditional expectation on the left-hand side of (5.40). By the scaling property of $\Phi$ (5.29), (5.30), we have that it can be written as

$$(5.43) \qquad \mathbb{E}[\Phi(1, \theta; Y^t) | \gamma] = \mathbb{E}[\Phi(1, \theta; Y) | \gamma^t].$$

(From now on we write $\mathbb{P}_\infty$ and $\mathbb{E}_\infty$ as $\mathbb{P}$ and $\mathbb{E}$, resp., using the subscript only for finite initial points.) For computing the right-hand side of (5.43), we first condition on $Y(1)$ and on whether or not there is a jump of $Y$ in $[1, 1 + \theta]$, and then if there is, at which time point it takes place. We get from that and (5.31), (5.32)

$$\mathbb{E}[\Phi(1, \theta; Y)|\gamma^t]$$

$$(5.44) \qquad = \Psi_1(\theta) \sum_{x \in \mathbb{N}^*} \mathbb{P}(Y_1 = x|\gamma^t) e^{-\theta t/\gamma_x}$$

$$(5.45) \qquad + \sum_{x \in \mathbb{N}^*} \mathbb{P}(Y_1 = x|\gamma^t) \int_0^\theta \frac{t}{\gamma_x} e^{-st/\gamma_x} \mathbb{E}[\Xi(x, s, Y_{[0, \theta-s]})|\gamma^t]\, ds,$$

where the sum can be taken in $\mathbb{N}^*$ due to Lemma 3.15, and we have used time homogeneity of $Y$. (We have made notation more compact by substituting parentheses with subscripts.)

We first note that the sum in (5.44) equals $\Lambda_t(\theta)$. Indeed

$$\Lambda_t(\theta) = \sum_{x \in \mathbb{N}^*} \mathbb{P}(Y_t = x|\gamma)\, \mathbb{P}_x(\text{no jump of } Y \text{ in } [t, t+\theta t]|\gamma)$$

$$(5.46)$$

$$= \sum_{x \in \mathbb{N}^*} \mathbb{P}(Y_t = x|\gamma)\, e^{-\theta t/\gamma_x} = \sum_{x \in \mathbb{N}^*} \mathbb{P}(Y_1 = x|\gamma^t)\, e^{-\theta t/\gamma_x},$$

where we have used the fact alluded to in Remark 5.19 above in the third equality. We now write the expression in (5.45) as

$$(5.47) \qquad \int_0^\theta \mathbb{E}[\Xi(1, s, Y_{[0, \theta-s]})|\gamma^t] \sum_{x \in \mathbb{N}^*} \mathbb{P}(Y_1 = x|\gamma^t) \frac{t}{\gamma_x} e^{-st/\gamma_x}\, ds$$

plus an error that is bounded above by

$$(5.48) \qquad \sup_{x \in \mathbb{N}^*, s \in (0, \theta)} |\mathbb{E}[\Xi(1, s, Y_{[0, \theta-s]})|\gamma^t] - \mathbb{E}[\Xi(x, s, Y_{[0, \theta-s]})|\gamma^t]|.$$

From (5.34), the absolute value of the difference of expectations in (5.48) can be bounded above by constant times

$$\sup_{x \in \mathbb{N}^*} \mathbb{P}[1, x \notin Y_{[0, \theta]}|\gamma^t] = \sup_{x \in \mathbb{N}^*} \mathbb{P}[1, x \notin Y_{[0, \theta t]}|\gamma]$$

$$(5.49)$$

$$= \sup_{x \in \mathbb{N}^*} \mathbb{P}[\Gamma^{(1, x)}(t\theta) < \sigma_1^1 \vee \sigma_1^x|\gamma],$$

where for $s > 0$

$$(5.50) \qquad \Gamma^{(1, x)}(s) = \sum_{y \neq 1, x} \gamma_x \sum_{i=1}^{N_s^{(y)}} T_i^{(y)}.$$

Thus, the right-hand side of (5.49) is bounded above by

$$(5.51) \qquad \mathbb{P}[\Gamma(t\theta) < T'|\gamma],$$

where $T'$ is a continuous random variable independent of $\Gamma$. It is clear from the fact that $\lim_{s \to 0} \Gamma(s) = 0$ that (5.51) vanishes as $t \to 0$.

We thus only have to consider (5.47). It can be written as

$$(5.52) \quad \int_0^\theta \Psi_2(s,\theta) \sum_{x \in \mathbb{N}^*} \mathbb{P}(Y_1 = x | \gamma^t) \frac{t}{\gamma_x} e^{-st/\gamma_x} \, ds = - \int_0^\theta \Psi_2(s,\theta) \Lambda_t'(s) \, ds,$$

where $\Lambda_t'(s) = \frac{d}{ds} \Lambda_t(s)$ (it is a straightforward exercise to show that the differentiation sign commutes with the sum), plus an error that is bounded above by

$$(5.53) \quad \begin{aligned} &\sup_{s \in (0,\theta)} |\mathbb{E}[\Xi(1,s,Y_{[0,\theta-s]})|\gamma^t] - \Psi_2(s,\theta)| \\ &= \sup_{s \in (0,\theta)} |\mathbb{E}[\Xi(1,s,Y_{[0,\theta-s]})|\gamma^t] - \Xi(1,s,\hat{\eta}_{\infty,\theta-s})|. \end{aligned}$$

From (5.34), (5.35), given $\varepsilon > 0$, there exists $\delta > 0$ such that the difference in (5.53) can be bounded above by constant times

$$(5.54) \quad \varepsilon + \mathbb{P}\left( \sup_{0 \le s \le \theta} \text{dist}(Y(s), \infty) > \delta \Big| \gamma^t \right).$$

Since, under $\gamma^t$, $Y_{[0,\theta]}$ converges in the sup norm to the identically in $[0,\theta]$ infinity path as $t \to 0$, we conclude that the expression in (5.53) vanishes as $t \to 0$.

We are thus left with taking the limit of

$$(5.55) \quad \begin{aligned} &\int_0^\theta \Psi_2(s,\theta) \Lambda_t'(s) \, ds \\ &= \Psi_2(\theta,\theta) \Lambda_t(\theta) - \Psi_2(0,\theta) - \int_0^\theta \Psi_2'(s,\theta) \Lambda_t(s) \, ds \end{aligned}$$

as $t \to 0$, where we have used the assumptions we made on $\Psi_2$ (5.36)–(5.37) to integrate by parts; note that $\Lambda_t(0) \equiv 1$. Collecting (5.44)–(5.55) and the above arguments together with the $L_1$ assumption (5.37) on $\Psi_2'$, the result then follows by (5.39) and dominated convergence, since $\Lambda_t$ is bounded (by 1). $\square$

REMARK 5.20. An alternative, longer argument for the validity of (5.40) in the form (5.41) for *almost every* $\gamma$, which has the advantage of requiring neither the uniformity assumption in (5.35) nor assumptions (5.36), (5.37)—see Remark 5.18 above—is to establish the convergence of $\Lambda_t'$ as $t \to 0$ to a (deterministic) function $\Lambda'$ (which turns out to be the derivative of $\Lambda$) for almost every $\gamma$. This can be done in an entirely similar fashion as in the proof of Theorem 5.11 below. We leave the details for the interested reader.

*Laplace transforms and Theorem* 5.11.   Before proving Theorem 5.11, we start with a simpler argument (at this point) for a weaker result, namely the a.s. convergence of a (double) Laplace transform of $\mathbb{E}[\Phi_1(\cdot,\cdot;Y)|\gamma]$. This requires the construction and results of Section 2 only. We shall then indicate how this computation can be used to identify the limit $\Lambda(\theta)$. Consider the function $c_\lambda(\mu)$ defined in (2.15). We can represent it as follows:

$$
\begin{aligned}
(5.56) \qquad c_\lambda(\mu) &= \lambda\mu \int_0^\infty \int_0^\infty e^{-\lambda s} e^{-\mu t} \mathbb{E}[\Phi_1(s,t;Y)|\gamma]\, ds\, dt \\
&= \lambda\mu \int_0^\infty \int_0^\infty e^{-\lambda s} e^{-\mu t} \Lambda_s\left(\frac{t}{s}\right) ds\, dt,
\end{aligned}
$$

with $c_\lambda$ defined in (2.15). For an aging result, it is natural to take $\lambda = \theta\mu$, and then take the limit as $\lambda \to \infty$. From (2.16), we have

$$
(5.57) \qquad c_\lambda(\lambda/\theta) = \frac{\sum_x (\lambda\gamma(x)/(1+\lambda\gamma(x))\lambda\gamma(x)/(\theta+\lambda\gamma(x)))}{\sum_x (\lambda\gamma(x)/(1+\lambda\gamma(x)))}.
$$

Before taking the limit, we note that both sums in (5.57) can be seen as sums over the increments of the Lévy process $V$ in $[0,1]$ [see paragraph of (5.8) above] of a function of the increments. We thus have by the scale invariance property of $V$ that the right-hand side of (5.57) for every $\lambda > 0$ has the same distribution as

$$
(5.58) \qquad \frac{\sum_{y\in[0,\lambda^\alpha]} (\gamma'(y)/(1+\gamma'(y)))(\gamma'(y)/(\theta+\gamma'(y)))}{\sum_{y\in[0,\lambda^\alpha]} \gamma'(y)/(1+\gamma'(y))},
$$

where the sum is over the increments $\{\gamma'_x\}$ of $V$ in $[0,\lambda^\alpha]$. Multiplying each factor of the quotient on the right-hand side of (5.58) by $\lambda^{-\alpha}$ and taking the limit as $\lambda \to \infty$, we get by the law of large numbers for sums of i.i.d. integrable variables that that quotient converges almost surely to

$$
\begin{aligned}
(5.59) \qquad & \frac{\mathbb{E}\sum_{y\in[0,1]} (\gamma'(y)/(1+\gamma'(y)))(\gamma'(y)/(\theta+\gamma'(y)))}{\mathbb{E}\sum_{y\in[0,1]} \gamma'(y)/(1+\gamma'(y))} \\
& = \frac{\int_0^\infty (w/(1+w))(w/(\theta+w))w^{-1-\alpha}\, dw}{\int_0^\infty (w/(1+w))w^{-1-\alpha}\, dw}
\end{aligned}
$$

as $\lambda \to \infty$.

REMARK 5.21.   That in principle says that the convergence of $c_\lambda(\lambda/\theta)$ as $\lambda \to \infty$ holds in probability; standard arguments relying on large deviation estimates for the sums on the right-hand side of (5.59) imply convergence almost everywhere.

In the forthcoming discussion, we shall use the notation

$$c(\theta) = \frac{\int_0^\infty (w/(1+w))(w/(\theta+w))w^{-1-\alpha}\,dw}{\int_0^\infty (w/(1+w))w^{-1-\alpha}\,dw}.$$

We have thus proved that, for a given value of $\theta$, $c_\lambda(\lambda/\theta)$ converges to $c(\theta)$. The monotonicity of $c_\lambda$ and the continuity of $c$ imply that we can in fact choose a set of $\gamma$'s of full measure such that the convergence holds for all $\theta > 0$ simultaneously.

Although the above computation does not seem to be sufficient to prove the convergence of $\Lambda_t(\theta)$, it can be used to identify its limit, as stated next.

PROPOSITION 5.22. *Suppose that, for (Lebesgue-almost) every $\theta \geq 0$, $\Lambda_t(\theta)$ converges to some limit $\Lambda(\theta)$ as $t \to 0$ (as stated in Theorem 5.11). Then $\Lambda$ equals $\Lambda^0$ given in (5.63) below.*

PROOF. After a simple change of variables in (5.56) we observe that

$$(5.60) \qquad c_\lambda(\lambda/\theta) = \int_0^\infty \int_0^\infty s e^{-s} e^{-st} \Lambda_{s/\lambda}(\theta t)\,ds\,dt.$$

The dominated convergence theorem then implies that

$$(5.61) \quad c_\lambda(\lambda/\theta) \to \int_0^\infty \int_0^\infty s e^{-s} e^{-st} \Lambda(\theta t)\,ds\,dt = \int_0^\infty \frac{1}{(1+t)^2}\Lambda(\theta t)\,dt.$$

We thus conclude that

$$(5.62) \qquad \begin{aligned} &\int_0^\infty \frac{1}{(1+t)^2}\Lambda(\theta t)\,dt \\ &= c(\theta) = \frac{\int_0^\infty (w/(1+w))(w/(\theta+w))w^{-1-\alpha}\,dw}{\int_0^\infty (w/(1+w))w^{-1-\alpha}\,dw}. \end{aligned}$$

The solution to (5.62) is easily seen to be given by the arcsine law. Indeed let

$$(5.63) \qquad \Lambda^0(\theta) = \frac{\sin(\pi\alpha)}{\pi}\int_{\theta/(1+\theta)}^1 s^{-\alpha}(1-s)^{\alpha-1}\,ds$$

be the distribution function of the arcsine law. Observe that

$$(5.64) \qquad \begin{aligned} &\int_0^\infty \frac{1}{(1+t)^2}\,dt \int_{\theta t/(1+\theta t)}^1 s^{-\alpha}(1-s)^{\alpha-1}\,ds \\ &= \int_0^\infty \frac{1}{(1+t)^2}\,dt \int_{\theta t}^\infty w^{-\alpha}\frac{dw}{1+w} \end{aligned}$$

$$(5.65) \qquad = \int_0^\infty \frac{w}{\theta+w}\,w^{-\alpha}\frac{dw}{1+w} \qquad \text{by Fubini,}$$

and therefore

$$(5.66) \qquad \int_0^\infty \frac{1}{(1+t)^2} \Lambda^0(\theta t)\, dt = c(\theta).$$

The next step is to invert (5.66), that is, prove that the transform

$$(5.67) \qquad \int_0^\infty \frac{1}{(1+t)^2} \phi(\theta t)\, dt$$

uniquely determines the probability distribution function $\phi$. This is an easy exercise in analysis: one can, for instance, write

$$(5.68) \qquad \frac{1}{(1+t)^2} = \int_0^\infty s\, e^{-s(1+t)}\, ds$$

and invert the two Laplace transforms.

We conclude that

$$(5.69) \qquad \Lambda(\theta) = \Lambda^0(\theta). \qquad\qquad \square$$

PROOF OF THEOREM 5.11. It is enough to get the result for a fixed $\theta > 0$. That we can find a full measure set of $\gamma$'s, such that the result holds for all $\theta > 0$ simultaneously, follows from the monotonicity of $\Phi_1$ and the continuity of $\Lambda$ in $\theta$.

We now give a full argument (independent of the above one). This argument uses the construction and results of Section 3 only.

The argument relies on an estimate for

$$(5.70) \qquad \mathbb{P}(Y_t = x | \gamma).$$

For $x \in \mathbb{N}^*$, $\{Y_t = x\}$ can be decomposed in the disjoint union of

$$(5.71) \qquad \{\Gamma^{(x)}(\sigma_1^{(x)}) \le t, \Gamma^{(x)}(\sigma_1^{(x)}) + \gamma_x T_1^{(x)} > t\},$$

where $\Gamma^{(x)} := \Gamma_0^{(x)}$ as in (3.34), and an event where $\gamma_x T_1^{(x)} \le t$ and $\Gamma^{(x)}(\sigma_2^{(x)}) \le t$. We thus have

$$(5.72) \qquad \begin{aligned} &|\mathbb{P}(Y_t = x|\gamma) - \mathbb{P}(\Gamma^{(x)}(S_1) \le t, \Gamma^{(x)}(\sigma_1^{(x)}) + \gamma_x T_1^{(x)} > t|\gamma)| \\ &\qquad \le (1 - e^{-t/\gamma_x})\mathbb{P}(\Gamma^{(x)}(S_2) \le t|\gamma), \end{aligned}$$

where $S_1$ and $S_2 - S_1$ are i.i.d. rate 1 exponentials which are independent of all other random variables around.

To establish the result we will prove the two following assertions:

$$(5.73) \qquad \sum_{x \in \mathbb{N}^*} e^{-\theta t/\gamma_x} \mathbb{P}(\Gamma^{(x)}(S_1) \le t, \Gamma^{(x)}(\sigma_1^{(x)}) + \gamma_x T_1^{(x)} > t|\gamma) \to \Lambda(\theta),$$

$$(5.74) \qquad \sum_{x \in \mathbb{N}^*} e^{-\theta t/\gamma_x} \mathbb{P}(\Gamma^{(x)}(S_2) \le t|\gamma) \to 0,$$

as $t \to 0$ for almost every $\gamma$. We rewrite the probability in (5.73) as follows:

$$(5.75) \qquad \mathbb{P}(\Gamma^{(x)}(S_1) \leq t | \gamma) - \mathbb{P}(\Gamma^{(x)}(\sigma_1^{(x)}) + \gamma_x \, T_1^{(x)} \leq t | \gamma),$$

and note that the second term equals

$$(5.76) \qquad \int_0^t e^{-(t-s)/\gamma_x} \mathbb{P}(\Gamma^{(x)}(S_1) \leq s | \gamma) \, ds$$
$$= \int_0^1 e^{-(1-s)t/\gamma_x} \mathbb{P}(\Gamma^{(x)}(S_1) \leq st | \gamma) \, ds.$$

Substituting in the left-hand side of (5.73), one sees that in order to prove the convergence in that display, it is enough to establish

$$(5.77) \qquad \sum_{x \in \mathbb{N}^*} e^{-\theta t/\gamma_x} \mathbb{P}(\Gamma^{(x)}(S_1) \leq t | \gamma) \to \hat{\Lambda}(\theta),$$

$$(5.78) \qquad \int_0^1 \sum_{x \in \mathbb{N}^*} e^{-((1+\theta)-s)t/\gamma_x} \mathbb{P}(\Gamma^{(x)}(S_1) \leq st | \gamma) \, ds \to \tilde{\Lambda}(\theta),$$

as $t \to 0$ for almost every $\gamma$, where $\hat{\Lambda}$ and $\tilde{\Lambda}$ are functions of $\theta$ only to be given below [see (5.89) and (5.94)]; we then have $\Lambda = \hat{\Lambda} - \tilde{\Lambda}$.

REMARK 5.23. We note that the left-hand sides of (5.77), (5.78) are both bounded above by $\sum_{x \in \mathbb{N}^*} e^{-\theta t/\gamma_x}$, which is almost surely finite for every $\theta, t > 0$, since $\sum_{x \in \mathbb{N}^*} \gamma_x$ is almost surely finite. They are thus almost surely finite.

We now observe that for almost every $\gamma$, $\Gamma^{(1)} \leq \Gamma^{(x)} \leq \Gamma$ for all $x \in \mathbb{N}^*$, where the first domination is a stochastic one (given $\gamma$), and follows from the decreasing monotonicity of $\gamma$.

To get (5.74), it suffices then to prove that for almost every $\gamma$

$$(5.79) \qquad \mathbb{P}(\Gamma^{(1)}(S_2) \leq t | \gamma) \sum_{x \in \mathbb{N}^*} e^{-\theta t/\gamma_x} \to 0 \qquad \text{as } t \to 0.$$

For (5.77), (5.78), it suffices to prove that for $i = 0, 1$

$$(5.80) \qquad \mathbb{P}(\Gamma^{(i)}(S_1) \leq t | \gamma) \sum_{x \in \mathbb{N}^*} e^{-\theta t/\gamma_x} \to \hat{\Lambda}(\theta),$$

$$(5.81) \qquad \int_0^1 \mathbb{P}(\Gamma^{(i)}(S_1) \leq st | \gamma) \sum_{x \in \mathbb{N}^*} \frac{t}{\gamma_x} e^{-((1+\theta)-s)t/\gamma_x} \, ds \to \tilde{\Lambda}(\theta),$$

as $t \to 0$, where $\Gamma^{(0)} = \Gamma$.

The next step is to replace, for $0 < s \leq 1$, $i = 0, 1$, $j = 1, 2$, $\mathbb{P}(\Gamma^{(i)}(S_j) \leq st | \gamma)$ by constant times $\omega_{ij}((st)^{-1})$, where for $r > 0$ $\omega_{ij}(r) := \mathbb{E}(\exp\{-r\,\Gamma^{(i)}(S_j)\} | \gamma)$.

This relies on a Tauberian theorem (see Theorem 3, Section 5, Chapter XIII of [14]), stating that as $t \to 0$, the quotient of the former quantity to the latter one converges to $1/\mathcal{G}(j\alpha)$ provided that for almost every $\gamma$

$$(5.82) \qquad \lim_{r \to \infty} \frac{\omega_{ij}(qr)}{\omega_{ij}(r)} = q^{-j\alpha} \qquad \text{for all } q > 0,$$

where, for $a > 0$, $\mathcal{G}(a) = \int_0^\infty t^a e^{-t}\, dt$.

Equation (5.82) is established in Lemma 5.29 below. From Lemma 5.27 and (5.99), we have that

$$(5.83) \qquad 0 \le \omega_{12}(t^{-1}) \sum_{x \in \mathbb{N}^*} e^{-\theta t/\gamma_x} \le \frac{\sum_x e^{-\theta/t^{-1}\gamma_x}}{(\sum_x t^{-1}\gamma_x/(1+t^{-1}\gamma_x))^2}.$$

Arguing as in the sentences above (5.58), we have that the right-hand side of (5.83) for every $t > 0$ has the same distribution as

$$(5.84) \qquad \frac{\sum_{y \in [0,t^{-\alpha}]} e^{-\theta/\gamma_y'}}{(\sum_{y \in [0,t^{-\alpha}]} \gamma_y'/(1+\gamma_y'))^2} = t^\alpha \frac{t^\alpha \sum_{y \in [0,t^{-\alpha}]} e^{-\theta/\gamma_y'}}{(t^\alpha \sum_{y \in [0,t^{-\alpha}]} \gamma_y'/(1+\gamma_y'))^2},$$

where the sum is over the increments $\{\gamma_x'\}$ of $V$ in $[0,t^{-\alpha}]$. Now the law of large numbers says that each factor in the quotient on the right-hand side of (5.84) converges almost surely to positive finite numbers as $t \to 0$. The extra factor of $t^\alpha$ in front of that expression then makes it vanish in that limit. That the same holds for the right-hand side of (5.83) follows as in Remark 5.21 above.

To get (5.80), we again need only get the limit for

$$(5.85) \qquad \omega_{i1}(t^{-1}) \sum_{x \in \mathbb{N}^*} e^{-\theta t/\gamma_x},$$

which by Lemma 5.27 and (5.99) is bounded above and below by

$$(5.86) \qquad \frac{\sum_x e^{-\theta/t^{-1}\gamma_x}}{k + \sum_x t^{-1}\gamma_x/(1+t^{-1}\gamma_x)},$$

$k = 0$ and $1$, respectively. As in (5.83)–(5.84), for every $t > 0$, (5.86) has the same distribution as

$$(5.87) \qquad \frac{\sum_{y \in [0,t^{-\alpha}]} e^{-\theta/\gamma_y'}}{k + \sum_{y \in [0,t^{-\alpha}]} \gamma_y'/(1+\gamma_y')}.$$

We can now multiply each term in the above quotient by $t^\alpha$ and take the limit as $t \to 0$. By the law of large numbers, this almost surely equals

$$(5.88) \qquad \frac{\mathbb{E}(\sum_{y \in [0,1]} e^{-\theta/\gamma_y'})}{\mathbb{E}(\sum_{y \in [0,1]} \gamma_y'/(1+\gamma_y'))} = \frac{\int_0^\infty e^{-\theta/w} w^{-(1+\alpha)}\, dw}{\int_0^\infty (w/(1+w)) w^{-(1+\alpha)}\, dw}.$$

An analogue of Remark 5.21 holds also here. We thus have from the above that

$$(5.89) \qquad \hat{\Lambda}(\theta) = \frac{1}{\mathcal{G}(\alpha)} \frac{\int_0^\infty e^{-\theta/w} w^{-(1+\alpha)} \, dw}{\int_0^\infty (w/(1+w)) w^{-(1+\alpha)} \, dw}.$$

It remains to get (5.81). Since $\mathbb{P}(\Gamma^{(i)}(S_j) \le st|\gamma)/\omega_{ij}((st)^{-1}) \to 1/\mathcal{G}(\alpha)$ as $t \to 0$ uniformly in $s \in (0,1]$, it suffices to get the limit for

$$(5.90) \qquad \int_0^1 \omega_{i1}((st)^{-1}) \sum_{x \in \mathbb{N}^*} \frac{t}{\gamma_x} e^{-((1+\theta)-s)t/\gamma_x} \, ds,$$

which by Lemma 5.27 and (5.99) reduces to getting the limit for

$$(5.91) \qquad \int_0^1 \frac{\sum_x (t/\gamma_x) e^{-((1+\theta)-s)t/\gamma_x}}{k + \sum_x (st)^{-1} \gamma_x/(1+(st)^{-1}\gamma_x)} \, ds,$$

$k = 0, 1$. It is clear that the quotient in the above integral, call it $\hat{\Lambda}_{s,t}(\theta)$, is bounded above by $\hat{\Lambda}_{1,t}(\theta)$, and the latter converges as $t \to 0$ almost surely [to $\hat{\Lambda}(\theta)$, as we just saw]. It thus suffices to establish the almost sure limit of $\hat{\Lambda}_{s,t}(\theta)$ as $t \to 0$ (independent of $k = 0, 1$) for every $s \in (0,1]$. That works similarly as above: we first multiply numerator and denominator of $\hat{\Lambda}_{s,t}(\theta)$ by $(st)^\alpha$, and get

$$(5.92) \qquad \hat{\Lambda}_{s,t}(\theta) = s^\alpha \frac{t^\alpha \sum_x (t/\gamma_x) e^{-((1+\theta)-s)t/\gamma_x}}{k(st)^\alpha + (st)^\alpha \sum_x (st)^{-1} \gamma_x/(1+(st)^{-1}\gamma_x)}.$$

Now the sums in the numerator and denominator of (5.92) are (marginally) equidistributed with $\sum_{y \in [0, t^{-\alpha}]} \frac{1}{\gamma'_y} e^{-((1+\theta)-s)/\gamma'_y}$ and $\sum_{y \in [0, (st)^{-\alpha}]} \frac{\gamma'_y}{1+\gamma'_y}$, respectively [by the very same argument made in between (5.57)–(5.58)]. By the law of large numbers, the quotient in (5.92), with the sums replaced by the respective ones just mentioned, converges almost surely as $t \to 0$ to

$$(5.93) \qquad \frac{\mathbb{E}(\sum_{y \in [0,1]} (1/\gamma'_y) e^{-((1+\theta)-s)/\gamma'_y})}{\mathbb{E}(\sum_{y \in [0,1]} \gamma'_y/(1+\gamma'_y))} = \frac{\int_0^\infty e^{-((1+\theta)-s)/w} w^{-(2+\alpha)} \, dw}{\int_0^\infty (w/(1+w)) w^{-(1+\alpha)} \, dw}.$$

We have an analogue of Remark 5.21 here as well, and thus can conclude that so does the original quotient. We thus have from the above that

$$(5.94) \qquad \bar{\Lambda}(\theta) = \frac{1}{\mathcal{G}(\alpha)} \frac{\int_0^1 \int_0^\infty s^\alpha e^{-((1+\theta)-s)/w} w^{-(2+\alpha)} \, dw \, ds}{\int_0^\infty (w/(1+w)) w^{-(1+\alpha)} \, dw}.$$

The result for fixed $\theta > 0$ is thus established with

$$(5.95) \qquad \Lambda = \hat{\Lambda} - \bar{\Lambda},$$

with $\hat{\Lambda}$, $\bar{\Lambda}$ given in (5.89), (5.94), respectively.  $\square$

REMARK 5.24. By (5.46), we see that

$$\Lambda_t(\theta) = \mathbb{E}(e^{-\theta t/\gamma_{Y_t}} \,|\, \gamma), \tag{5.96}$$

so Theorem 5.11 and Proposition 5.22 establish that, for almost every $\gamma$, $t/\gamma_{Y_t}$ converges in distribution as $t \to 0$ to the random variable $Z$ whose Laplace transform $\mathbb{E}(e^{-\theta Z})$ is $\Lambda^0(\theta)$. This is another way of understanding the basic mechanism for the aging phenomenon in this process (there is no change for a time of order $t$ when the process has aged $t$ units of time). It is also a *macroscopic* version of the last assertion of Proposition 2.10 of [9].

REMARK 5.25. Lemma 2.11 of [9] establishes the continuity of the distribution of the random variable $Z$ in Remark 5.24. From (5.89), (5.94), one readily finds its density with respect to Lebesgue measure, given by

$$\frac{1}{\mathcal{G}(\alpha) \int_0^\infty w^{-\alpha} (1+w)^{-1} \, dw} z^{\alpha-1} \int_0^1 \alpha s^{\alpha-1} e^{-(1-s)z} \, ds, \qquad z > 0. \tag{5.97}$$

REMARK 5.26. The convergence in Remark 5.24 suggests that a stronger result holds, namely, that the process $\{Z_t^{(\varepsilon)} := \varepsilon^{-1} \gamma_{Y_{\varepsilon t}}, \, t \geq 0\}$, with $Y(0) = \infty$, and thus $Z_0^{(\varepsilon)} \equiv 0$, converges in distribution as $\varepsilon \to 0$ to a nontrivial limit for almost every $\gamma$. This result would contain Theorem 5.11. A related scaling limit, namely that of $\{\varepsilon^{-1} \Gamma(\varepsilon^\alpha t), \, t \geq 0\}$ as $\varepsilon \to 0$ (to an $\alpha$-stable subordinator), would also imply Theorem 5.11, and would be a *macroscopic* version of the corresponding *mesoscopic* convergence of the *clock process* proved in [3]. We believe both these limits can be established, for example, by checking standard convergence criteria; these are points of our current research.

LEMMA 5.27. *For $i = 0, 1$, $j = 1, 2$, and almost every $\gamma$,*

$$\omega_{ij}(r) = \left(1 + \sum_{x \neq i} \frac{r\gamma_x}{1 + r\gamma_x}\right)^{-j}. \tag{5.98}$$

PROOF. Exercise.  □

REMARK 5.28. The condition $x \neq 0$ in the sum in (5.98) (in the case when $i = 0$) is empty since $x \geq 1$. From (5.98), we have that for $i = 0, 1$, $j = 1, 2$

$$\left(1 + \sum_x \frac{r\gamma_x}{1 + r\gamma_x}\right)^{-j} \leq \omega_{ij}(r) \leq \left(\sum_x \frac{r\gamma_x}{1 + r\gamma_x}\right)^{-j}. \tag{5.99}$$

LEMMA 5.29. *Equation (5.82) holds for almost every $\gamma$.*

PROOF.  For a fixed $\lambda > 0$, and then for all rational $\lambda > 0$, it follows from (5.98) and a law of large number argument as in the proof of Theorem 5.11. The result for all $\lambda > 0$ can be argued from that, using the monotonicity of $\omega_{ij}(\cdot)$ and the continuity of the limit.  □

**Acknowledgments.**  We would like to thank Gérard Ben Arous for letting us have a preliminary version of [3] and Véronique Gayrard for many discussions and comments. We also thank an anonymous referee for many stimulating comments that helped improve the presentation and suggested interesting related points for further investigation (like the one about the latter convergence in Remark 5.26).

The authors benefitted from travel grants by the USP-COFECUB and the Brazil–France agreements.

IME-USP
RUA DO MATÃO 1010
05508-090 SÃO PAULO SP
BRAZIL
E-MAIL: lrenato@ime.usp.br

CMI
UNIVERSITÉ DE PROVENCE
39 RUE JOLIOT-CURIE
13013 MARSEILLE
FRANCE
E-MAIL: pierre.mathieu@cmi.univ-mrs.fr